\documentclass[10pt]{siamltex}

\usepackage{amssymb}
\usepackage{graphicx}
\usepackage{epstopdf}
\usepackage{caption2}
\usepackage{psfrag}
\usepackage{amsmath}
\usepackage{amscd}
\usepackage{amsfonts}
\usepackage{bm}   %
\usepackage{float}
\usepackage{latexsym}
\usepackage{lscape}
\usepackage{mathrsfs}
\usepackage{verbatim}
\usepackage{appendix}
\usepackage{tikz}
\allowdisplaybreaks[4]
\newtheorem{remark}{Remark}[section]
\newtheorem{thm}{Theorem}
\newtheorem{lem}[thm]{Lemma}
\newtheorem{exam}{Example}
\newcommand{\ds}{\displaystyle}
\newcommand{\f}{\frac}

\newcommand{\ba}{\begin{array}}
\newcommand{\ea}{\end{array}}
\newcommand{\beq}{\begin{equation}}
\newcommand{\eeq}{\end{equation}}
\newcommand{\bex}{\begin{eqnarray*}}
\newcommand{\eex}{\end{eqnarray*}}
\newcommand{\be}{\begin{equation*}}
\newcommand{\ee}{\end{equation*}}

\def\br{\begin{eqnarray}}
\def\er{\end{eqnarray}}
\def\brr{\bq\begin{array}{rlll}}
\def\err{\end{array}\eq}

\def\N{{\mathbb {N} }}

\def\d{{\mathrm {d} }}

 \makeatletter
\def\ps@pprintTitle{%
 \let\@oddhead\@empty
 \let\@evenhead\@empty
 \def\@oddfoot{\centerline{\thepage}}%
 \let\@evenfoot\@oddfoot}
\makeatother
\usepackage[mathlines]{lineno}

\title{A fast solver for spectral element approximation applied to fractional differential equations using hierarchical matrix approximation\thanks{This work was supported by the MURI/ARO on ``Fractional PDEs for Conservation Laws and Beyond: Theory, Numerics and Applications" (W911NF-15-1-0562).} The research of the first author was partially supported by NSF of China
(11771083, 61672005), NSF of Fujian Province (2016J01013).
}

\author{
  Xianjuan Li\footnotemark[2]
  \and
  Zhiping Mao\footnotemark[3] \footnotemark[5]
  \and
  Fangying Song\footnotemark[3]
  \and
  Hong Wang\footnotemark[4]
  \and
  George Em Karniadakis\footnotemark[3]
}

\begin{document}

\markboth{}{}

\footnotetext[2]{College of Mathematics and Computer Science, Fuzhou University, Fuzhou, 350108 China ({xjli@fzu.edu.cn}). }
\footnotetext[3]{Division of Applied Mathematics, Brown University, Providence, RI 02912, USA ({zhiping\_mao@brown.edu}, {fangying\_song@brown.edu}, {george\_karniadakis@brown.edu}).}
\footnotetext[4]{Department of Mathematics, University of South Carolina, Columbia, SC 29208, USA ({hwang@math.sc.edu}).}
\footnotetext[5]{Corresponding author.}

\maketitle

\begin{abstract}
We develop a fast solver for the spectral element method (SEM) applied to the two-sided fractional diffusion equation on uniform, geometric and graded meshes.
By approximating the singular kernel with a degenerate kernel, we construct a hierarchical matrix (H-matrix) to represent the stiffness matrix of the SEM and provide error estimates verified numerically.
We can solve efficiently  the H-matrix approximation problem using a hierarchical LU decomposition method, which reduces the computational cost to $O(R^2 N_d  \log^2N) +O(R^3 N_d  \log N)$, where $R$ it is the rank of submatrices of the H-matrix approximation,  $N_d$ is the total number of degrees of freedom and $N$ is the number of elements.
However, we lose the high accuracy of the SEM. Thus, we solve the corresponding preconditioned system by using the H-matrix approximation problem as a preconditioner, recovering the high order accuracy of the SEM. The condition number of the preconditioned system is independent of the polynomial degree $P$ and grows with the number of elements, but at modest values of the rank $R$ is below order 10 in our experiments, which represents a reduction of more than 11 orders of magnitude from the unpreconditioned system; this reduction is higher in the two-sided fractional derivative compared to one-sided fractional derivative. The corresponding cost is $O(R^2 N_d \log^2 N)+O(R^3 N_d  \log N)+O(N_d^2)$. Moreover, by using a structured mesh (uniform or geometric mesh), we can further reduce the computational cost to  $O(R^2 N_d\log^2 N) +O(R^3 N_d  \log N)+ O(P^2 N\log N)$ for the preconditioned system.
We present several numerical tests to illustrate the proposed algorithm using $h$ and $p$ refinements.
\end{abstract}

\begin{keywords}
Hierarchical Matrix, Spectral element method, Hierarchical LU decomposition, Preconditioner, Toeplitz matrix
\end{keywords}

\begin{AMS}
  65N35, 65M70,  41A05, 41A25
\end{AMS}

\section{Introduction}
In this paper, we consider the following one-dimensional two-sided Riemann-Liouville fractional differential equation (FDE):
\begin{equation}\label{Model}
\begin{aligned}
\ds - D^2 \bigl (\theta ~{}_a I_x^{2-\alpha} u + (1-\theta) ~{}_xI_b^{2-\alpha} u \bigr ) + \rho u  = f(x), & \quad  x \in \Lambda, \\
u(a) = c_1,\; u(b) = c_2, &
\end{aligned}
\end{equation}
where $\Lambda : = (a,b)$, $1 < \alpha < 2$, $f$ is a given function and $c_1,\,c_2$ are two constants.
Here $Du(x) := u'(x)$ is the first-order differential operator, $0 \le \theta \le 1$ indicating the relative weight of forward versus backward transition probability. The left and right fractional integrals of order $0 < \beta < 1$ are defined by~\cite{podlubny1998fractional}
\begin{equation*}
    {}^{}_aI_x^{\beta}u(x) \ds := \f1{\Gamma(\beta)} \int_a^x \f{u(s)}{(x-s)^{1-\beta}}ds, \quad
{}^{}_xI^{\beta}_b u(x) \ds := \f1{\Gamma(\beta)} \int_x^b \f{u(s)}{(s-x)^{1-\beta}}ds,
\end{equation*}
where $\Gamma(\cdot)$ is the Gamma function.

The two-sided fractional diffusion equation is required in applications such as hydrology~\cite{benson2000application, chakraborty2009parameter, zhang2016bounded} and turbulent transport in plasmas~\cite{del-Castillo-Negrete2006}. More applications can be found in \cite{deng2006parameter, baeumer2007fractional, schumer2009fractional}.

Due to the fact that analytic forms of the solutions of  FDEs are difficult to obtain, numerical methods are favored to study the behaviour of anomalous diffusion.
In the past years, extensive research has been carried out on the development of numerical methods for FDEs, including finite difference methods~\cite{MeerMM04, DengweihuaMC15}, finite element methods (FEMs)~\cite{ervin2006variational, Wanghong13, zhao2017adaptive}, spectral methods~\cite{LiXu10, MohsenKar13, MaoShen16, MaoChenShen16,  ErvinRoop18, MaoKar18}, spectral elements methods (SEM)~\cite{Zay.K14, MaoShen17Adv} and references therein.
It is now well-known that the solutions of FDEs exhibit end-points singularities, which causes difficulties in developing high accuracy numerical methods. To obtain high accuracy solution, Jin et al. proposed a FEM by using a regularity reconstruction to improve the convergence rate~\cite{J.Z.2015ESIAM}. Zhao et al. developed an adaptive FEM to improve the accuracy of the numerical solutions~\cite{zhao2017adaptive}. However, the proposed FEMs are low order method using piecewise linear approximation. On the other hand, efficient and high order Petrov-Galerkin spectral (PGS) methods have been developed by using the Jacobi poly-fractonomials as basis functions~\cite{MohsenKar13, ChenShenWang2016MathComp, MaoChenShen16, MaoKar18}. For a single model problem (without reaction term), the resulted linear system is sparse (diagonal), and for homogeneous boundary conditions, the error only depends on the regularity of the right hand function $f$. Unfortunately, if a reaction term is presented or if we use non-homogeneous boundary conditions, these advantages do not hold anymore because the singular behavior at the end points is complex. To overcome these issues, a more flexible SEM was proposed by Mao and Shen showing that the convergence rate with geometric meshes is exponential with respect to the square root of the number of degrees of freedom without prior knowledge about the singular behavior~\cite{MaoShen17Adv}.

\begin{figure}[!t]
\centering
\includegraphics[width=0.48\textwidth,height=0.4\textwidth]{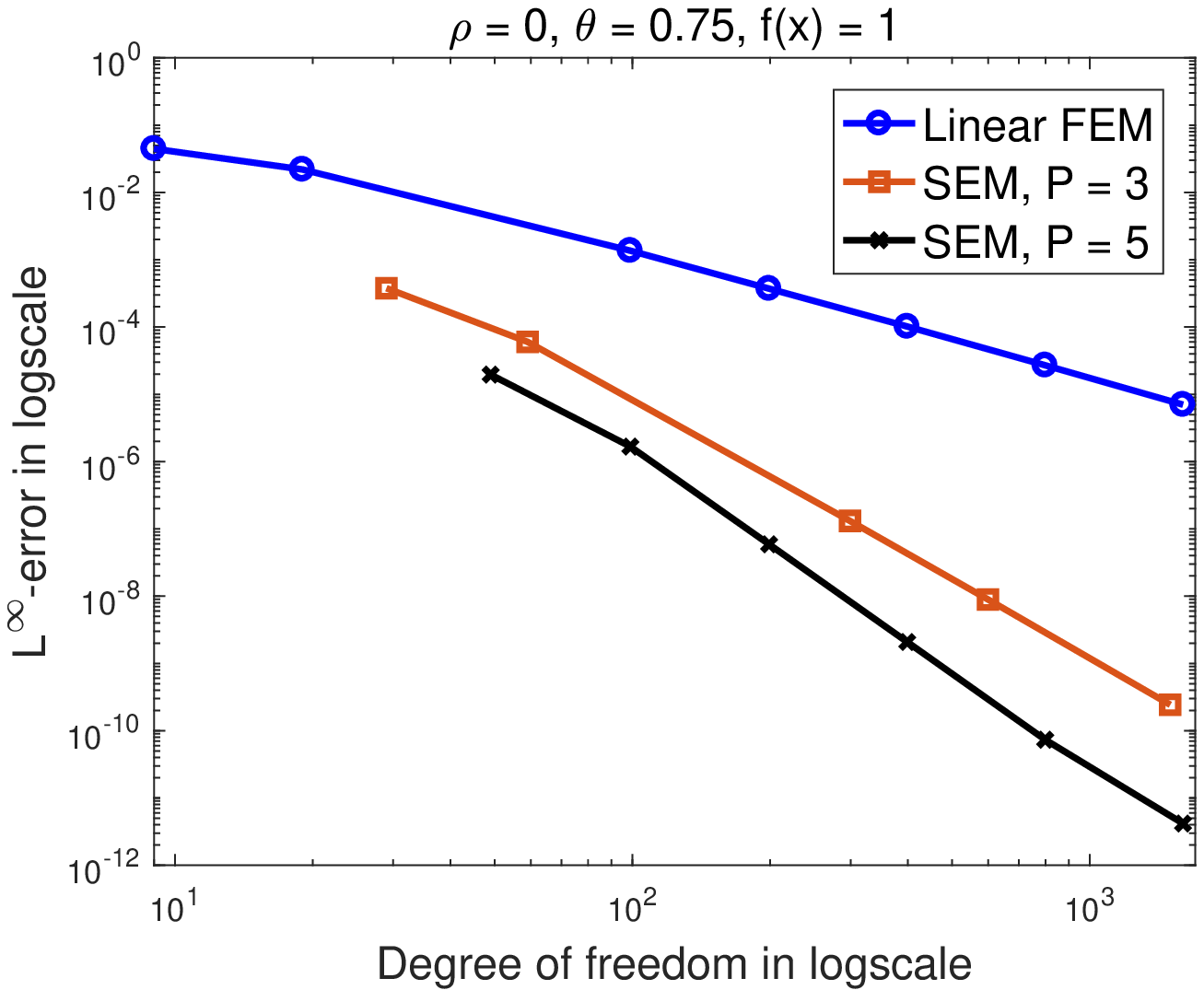}
\includegraphics[width=0.48\textwidth,height=0.4\textwidth]{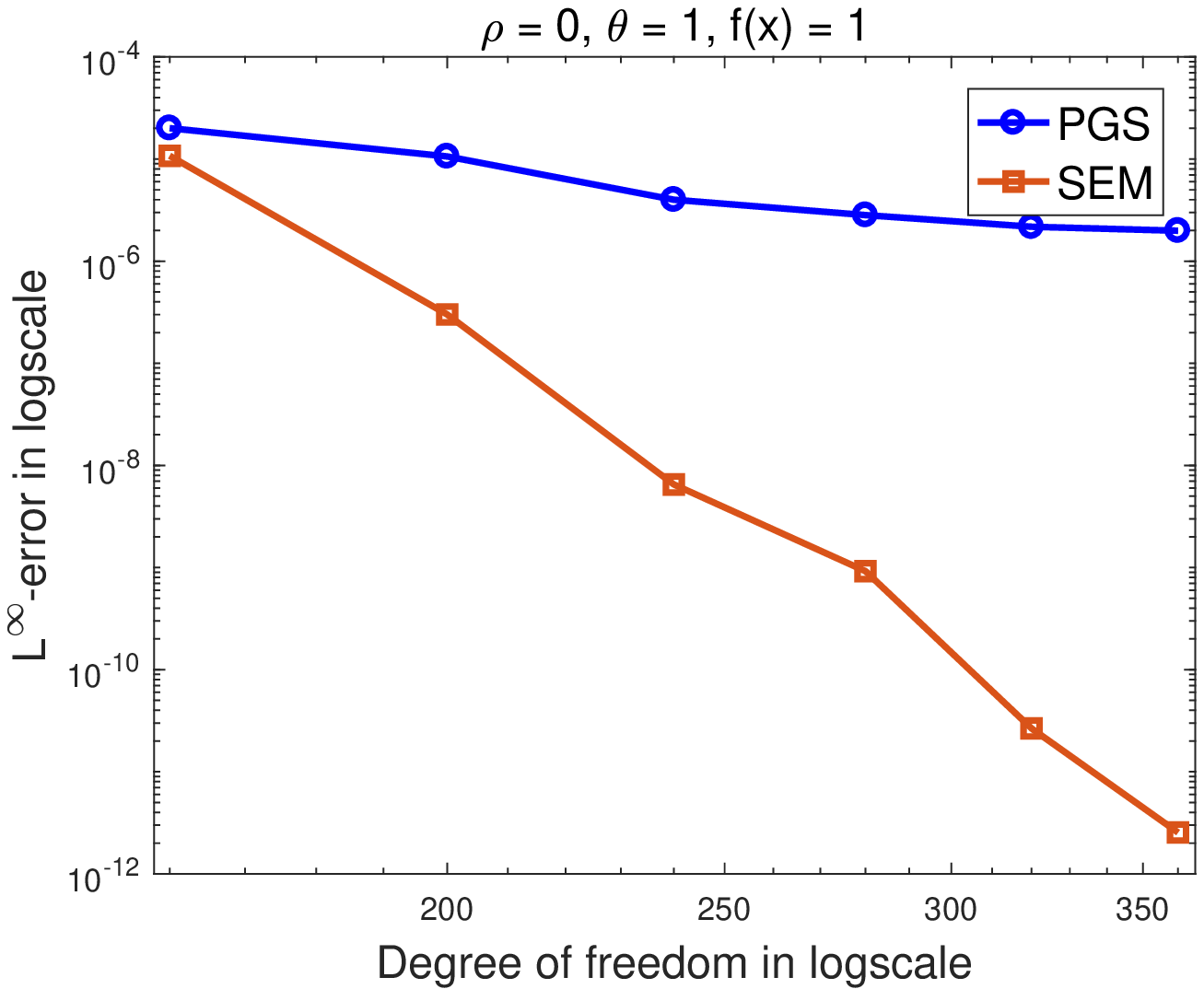}
\caption{\scriptsize $L^\infty$-error for $\alpha = 1.6$. $f(x) = 1$, $a = -1, b = 1$. Left: Comparison of the $L^\infty$-error obtained by the SEM against the linear FEM, $c_1 = c_2 = 0, \, \rho = 0,\, \theta = 0.75$. The mesh for both FEM and SEM is the same graded mesh with $x_i = (2i/N)^7 -1,i= 0,1,\ldots, N/2,\, x_{i} = 1-((N-2i)/N)^5,i= N/2+1,\ldots,N$. Right: Comparison of the $L^\infty$-error obtained by the SEM against the PGS, $c_1 = 0,\, c_2 = 3*2^{\alpha-1}/\gamma(1+\alpha), \, \rho = 0, \, \theta = 1$. The mesh for SEM is the graded mesh with $x_0 = -1, x_i =x_0+ 2(i/N)^5,i= 1,\ldots, N,\, N = 40$, and the polynomial degree varies from 4 to 9.
}
\label{fig:motivation} 
\end{figure}

In the present work, we adopt the SEM proposed in \cite{MaoShen17Adv}. To illustrate the high accuracy of the SEM, we show in Figure \ref{fig:motivation} the comparisons of $L^\infty$-error obtained by the SEM against a linear FEM with homogeneous boundaries (left) and a PGS developed in \cite{MaoKar18} with non-homogenous boundary conditions (right) for smooth-right-hand side, leading to a singular solution. Observe that the SEM can deliver much higher accuracy than the linear FEM and the PGS method.
However, the resulted linear system of the SEM is dense due to the non-locality of the fractional operators and the condition number may grow fast using a non-uniform mesh. Hence, we address here how to efficiently solve the dense linear system.
Clearly, it is too expensive to use the direct method since it requires the storage to be $O(N_d^2)$ and the computational cost to be $O(N_d^3)$, where $N_d$ is the number of the degrees of freedom.
For a uniform mesh, the discretization matrices are Toeplitz-like, and therefore the memory can be reduced to  $O(N)$ and the computational cost to $O(N\log N)$, where $N$ is the number of elements, using fast matrix-vector multiplication in conjunction with effective preconditioners, see \cite{wang2013superfast} and references therein.
However, for a general mesh, this kind of approach fails. To this end, Zhao et al.~\cite{zhao2017adaptive} developed a fast solver for a adaptive FEM with general mesh using hierarchical matrix approximation. Subsequently, Ainsworth and Glusa~\cite{AG2017CMAME} developed a fast solver for an adaptive FEM on a two-dimensional unstructured mesh with assembly, matrix-vector product and computation of the error estimators scaling quasi-linearly with respect to the number of unknowns.
To the best of our knowledge, the present work represents the first attempt to develop a fast solver for the SEM discretization on general one-dimensional meshes using hierarchical matrix (H-matrix) approximation.

The main contributions of this work are as follows:
\begin{itemize}
  \item We construct a H-matrix approximation for the high order SEM by replacing the singular kernel with a degenerate kernel and provide  the error estimates verified numerically for the approximation of entries of the stiffness matrix.
  \item We efficiently solve the H-matrix approximation problem by using hierarchical LU decomposition to reduce the storage to $O(R N_d\log N)$ and the computational cost to $O(R^2 N_d  \log^2N) +O(R^3 N_d  \log N)$, where $R$ is the rank of  submatrices of the H-matrix approximation.
  \item We employ the H-matrix approximation as a preconditioner to retain the high accuracy of the SEM and efficiently solve the preconditioned system at the cost of $O(R^2 N_d \log^2 N)+O(R^3 N_d  \log N)+O(N_d^2)$. We can further reduce the computational cost to $O(R^2 N_d\log^2 N)+O(R^3 N_d  \log N) + O(P^2 N\log N)$ if using structured (uniform or geometric) mesh, where $P$ is the polynomial  degree of each spectral element.
\end{itemize}

The  rest of the paper is organized as follows.
In section \ref{sec:wksem}, we present the weak formulation of \eqref{Model} and its spectral element discretization.
In section \ref{sec:HMR}, we construct the low rank H-matrix representation for the elements of the stiffness matrix.
Then we develop a fast solver for the resulted linear system in section \ref{sec:Fsolver}.
In section \ref{sec_num}, we illustrate the proposed algorithms by using numerical simulations.
Finally, we conclude in section \ref{sec:conclusion}.

\section{Weak formulation and spectral element approximation}\label{sec:wksem}

In this section we present the weak formulation and its spectral element approximation.

\subsection{Weak formulation}

Let $L_2(\Lambda)$ be the Hilbert space of Lebesgue square integrable functions on $\Lambda$. Let $H^\mu(\Lambda)$ ($\mu \ge 0$) be the fractional Sobolev space of order $\mu$ and $H^\mu_0(\Lambda)$, with $\mu > 1/2$, be its subspace that consists of functions with zero boundary conditions. Let $H^{-\mu}(\Lambda)$ be the dual space of $H^\mu_0(a,b)$. Ervin and Roop \cite{ervin2006variational} derived a Galerkin weak formulation for problem (\ref{Model}): Given $f \in H^{-\alpha/2}(\Lambda)$, find $u \in H^{\alpha/2}_0(\Lambda)$ such that
\begin{equation}\label{WeakForm}
a(u,v) = \langle f, v \rangle~~\forall v \in H^{\alpha/2}_0(\Lambda),
\end{equation}
where the bilinear form is given by
\begin{equation*}
    a(u,v) := \theta\left({}_aD_x^{\alpha/2} u, {}_xD_b^{\alpha/2} v \right) + (1-\theta)\left({}_xD_b^{\alpha/2} u, {}_aD_x^{\alpha/2} v \right) + (\rho u,v),
\end{equation*}
$(\cdot,\cdot)$ and $\langle \cdot, \cdot \rangle$ denote the inner product in the space $L_2(\Lambda)$ and the duality pairing of $H^{-\alpha/2}(\Lambda)$ and $H^{\alpha/2}_0(\Lambda)$. For the weak problem \eqref{WeakForm}, we have the following theorem~\cite{ervin2006variational, LiXu10}:
\begin{thm}\label{thm:Wellposedness}
The bilinear form $a(\cdot,\cdot)$ is coercive and continuous on the product space $H^{\alpha/2}_0(\Lambda) \times H^{\alpha/2}_0(\Lambda)$, i.e., there are positive constants $0 < C_1 < C_2 < \infty$ such that
\begin{equation}\label{Coercivity}
 a(w,w) \ge C_1 \|w \|_{H^{\alpha/2}}^2, \; 
| a(w,v) | \le C_2  \|w \|_{H^{\alpha/2}} \| v \|_{H^{\alpha/2}} \quad  \forall w, v \in H^{\alpha/2}_0(\Lambda).
\end{equation}
Hence, the Galerkin weak formulation (\ref{WeakForm}) has a unique solution $u \in H^{\alpha/2}_0(\Lambda)$ with the stability estimate
\begin{equation*}
    \ds \| u \|_{H^{\alpha/2}(\Lambda)} \le K \| f \|_{H^{-\alpha/2}(\Lambda)}.
\end{equation*}
\end{thm}

\subsection{Spectral element discretization}

We are now in the position to introduce a SEM for the model problem (\ref{Model}) based on the weak formulation (\ref{WeakForm}). In particular, we adopt the SEM proposed by Mao and Shen~\cite{MaoShen16}. We begin by defining a spatial partition on $\Lambda$
\begin{equation}\label{Partition}
a =: x_0 < x_1 < \ldots < x_i < \ldots < x_{N-1} < x_N := b,   \quad h_{i}=x_i-x_{i-1}.
\end{equation}
Then, we give a complete set of basis functions.
We first define the continuous and piecewise-linear \emph{nodal basis functions} $\{\varphi_{j}\}_{j=1}^{N-1}$ on $\Lambda$ with respect to the partition \eqref{Partition} such that $\varphi_{j}(x_k) = \delta_{jk}$, which equal to 1 if $j=k$ and 0 otherwise.
Then we construct the \emph{model basis functions}.
For any nonnegative integer $p$, let $L_p(x)$ be the $p$-th degree Legendre polynomial on the interval $[-1,1]$, which can be defined by the following recurrence relation \cite{karniadakis2013spectral,shen2011spectral}
\begin{equation*}
    L_0(x) =1,\; L_1(x)=x,\; L_{p+1}(x) =\ds\frac{2p+1}{p+1}xL_p(x)-\frac{p}{p+1}L_{p-1}(x),\quad p \geqslant 1.
\end{equation*}
It is well known that Legendre polynomials are orthogonal on the interval $[-1,1]$
$$\ds \int_{-1}^1 L_p(x) L_q(x) dx = \f{2}{2p+1}\delta_{pq}, \quad \ds L_p(\pm 1)  = (\pm 1)^p.$$
Hence, it is clear that the functions
\begin{equation}\label{RefBasis}
\psi_p(x) := L_{p-1}(x) - L_{p+1}(x), \, p\ge 1,
\end{equation}
are linearly independent satisfying $\psi_p(\pm 1) = 0$. With the transformation
$$ \ds \xi_i = \xi_i(x) := \f{2x- (x_{i-1}+x_i)}{x_i - x_{i-1}},$$
We define the following \emph{model basis functions} on $\Lambda$ supported  on $\Lambda_i: = [x_{i-1},x_i]$
\begin{equation}\label{LocalBasis}
\phi_{i}^{p}(x) := \bigg \{ \begin{array}{cc}
\ds \psi_p(\xi_i(x)), & x \in [x_{i-1},x_i], \\[0.05in]
0, & x \in \Lambda\backslash [x_{i-1},x_i].
\end{array}\end{equation}

Let $S_{N}^{P}(\Lambda)$ be the spectral element space defined by
\begin{equation*}
S_{N}^{P}(\Lambda) := \mathrm{span}\{\varphi_{n}\}_{n=1}^{N-1} \oplus  \mathrm{span}\{\phi_{1}^{p}\}_{p=1}^{P_1-1} \oplus \cdots \oplus \mathrm{span}\{\phi_{N}^{p}\}_{p=1}^{P_N-1}.
\end{equation*}
Then the SEM to \eqref{Model} can be formulated as follows: Find $u_{N}^{P} \in S_{N}^{P}(\Lambda)$ such that
\begin{equation}\label{SEM}
a(u_{N}^{P},v_{N}^{P}) = \langle f,v_{N}^{P} \rangle, \quad v_{N}^{P} \in S_{N}^{P}(\Lambda).
\end{equation}
Since $S_{N}^{P} \subset H^{\alpha/2}_0(\Lambda)$, Theorem \ref{thm:Wellposedness} ensures that the SEM \eqref{SEM} has a unique solution $u_{N}^{P} \in S_{N}^{P}(\Lambda)$.

Expanding $u_{N}^{P}(x)$ by
\begin{equation*}
    u_{N}^{P}(x) = \sum_{n=1}^{N-1} \tilde{u}_n \varphi_n(x) + \sum_{n=1}^N\sum_{p=1}^{P_n -1} \hat{u}_n^p \phi_n^p(x)
\end{equation*}
and using the standard procedure, we can obtain the following linear system
\beq\label{lsys}
A X := (\rho M+S)X =G,
\eeq
where $M$ is the mass matrix and
\begin{equation*}
    S = \theta S_l + (1-\theta)S_l^T,\text{ with }
    S_l = \left[
\ba{cc}
{B} & {F}\\
{E} & {C}
\ea
\right]
\end{equation*}
is the stiffness matrix, see Figure \ref{fig200}  (left) for the partition of $S_l$. Here, for $ 1\le i, j\le N,\, 1\le p\le P_i -1,\, 1\le q\le P_j -1$,
\begin{align*}
%
&B_{ij}^{pq} = \left({}_aD_x^{\alpha/2} \phi_{j}^q(x), {}_xD_b^{\alpha/2} \phi_{i}^p(x) \right), \quad
C_{ij} = \left({}_aD_x^{\alpha/2} \varphi_{j}(x), {}_xD_b^{\alpha/2}\varphi_{i}(x) \right),\\
&E_{ij}^{q} = \left({}_aD_x^{\alpha/2}\phi_{j}^q(x), {}_xD_b^{\alpha/2} \varphi_{i}(x) \right), \quad
F_{ij}^{p} = \left({}_aD_x^{\alpha/2}\varphi_{j}(x), {}_xD_b^{\alpha/2} \phi_{i}^p(x) \right)
\end{align*}
and
\begin{equation*}
    \begin{aligned}
     & X = \big[\hat{u}_{1}^{1}, \cdots, \hat{u}_{1}^{P_1-1}; \cdots ; \hat{u}_{N}^{1}, \cdots, \hat{u}_{N}^{P_N-1}; \tilde{u}_1, \cdots, \tilde{u}_{N-1}\big]^T,\; G = [\hat{G}, \tilde{G}]^T,\\
    & \hat{G} = \big[\hat{g}_{11},\cdots, \hat{g}_{1,P_1-1};\cdots;\hat{g}_{N1},\cdots, \hat{g}_{N,P_N-2};\big], \; \hat{g}_{jk} = \big(g,\phi_j^k(x)\big), \\
    & \tilde{G}=[\tilde{g}_1,\cdots,\tilde{g}_{N-1}], \; \tilde{g}_{j} = \big(g,\varphi_j(x)\big).
    \end{aligned}
\end{equation*}
The total number of degrees of freedom is $N_d:=\sum\limits_{i=1}^{N}P_i-1$.

\section{Hierarchical matrix representation}\label{sec:HMR}
It is known that the mass matrix $M$ is a sparse matrix and the stiffness matrix $S$ is a dense matrix  due to the non-locality of the fractional operators.
In this section, we construct a H-matrix, which  can be stored in a data-sparse format, to approximate the stiffness matrix $S$.
\subsection{Approximation of the kernel}
The main idea of H-matrix representation is to approximate the kernel
$k(t,s)= |t-s|^{1-\alpha}$ 
by a degenerate kernel $\tilde{k}(t,s)$, i.e.,
\beq\label{eqdegen}
\tilde{k}(t,s) = \sum_{\nu=0}^{R-1}  g_\nu(s) h_\nu(t).
\eeq
In order to approximate the stiffness matrix by low rank H-matrix,
there are two key properties: firstly, the approximation is degenerate, i.e.,
the variables $t$ and $s$ have to be separated, secondly, it has to converge rapidly to the original kernel.
However, due to the singularity of the kernel function  $k(t,s)= |t-s|^{1-\alpha}$  at the line $t=s$,
the degenerate function \eqref{eqdegen} converges very slowly  in the whole domain $[a, b] \times [a,b]$.
Thus, the local approximations  are made  for subdomains of $[a, b] \times [a,b]$  in which $k(t,s)$ is smooth.
The following admissibility condition is used to characterize this local property of the  approximation.
\begin{definition}\cite{borm2005hierarchical,hackbusch2015hierarchical}
Let $\tau=[a',b']$  and $\sigma=[c',d']$ be two bounded intervals,  $\tau  \times\sigma$ are admissible
if  the following  admissibility condition holds,
\beq\label{admis}
\max\{Diam(\sigma), Diam(\tau)\} \leq \lambda\cdot Dist(\tau, \sigma),  \quad  \lambda>0,
\eeq
 otherwise $\tau  \times\sigma$ are inadmissible.
\end{definition}

Obviously, the kernel $k(t,s)$ defined on admissible intervals is nonsingular.
Concretely,  the usual truncated  Taylor expansion matches the properties of degenerate function \eqref{eqdegen}.
This means that if $[a',b']  \times [c',d']$ are admissible,
then  for  $t\times s \in  [a',b']  \times [c',d']$, the truncated  Taylor expansion can be used
to approximate the $k(t,s)= |t-s|^{1-\alpha}$.
The following lemma  provides the error of $|k(t,s)-\tilde{k}(t,s)|$ for $t>s$ (in this case, $k(t,s) = (t-s)^{1-\alpha}$).
\begin{lem}
(Taylor series of $(t-s)^{1-\alpha}$)  Let $\tau= [a',b']$ and $  \sigma=[c',d']$  be be two bounded intervals such that $d'< a'$.
Assuming that $ \tau\times \sigma$ satisfy the admissibility condition \eqref{admis}, then for $t\times  s \in  \tau\times \sigma, \; R\in \N$, we can approximate the kernel $k(t,s) = (t-s)^{1-\alpha}$ by
\begin{equation}\label{eqn:Lexpan}
    \tilde{k}(t,s) = \sum_{\nu=0}^{R-1}  \frac{1}{\nu!} \partial_s^\nu [(t-s_0)^{1-\alpha}] (s-s_0)^{\nu}
    =\sum_{\nu=0}^{R-1}  g_\nu(s,s_0) h_\nu(t,s_0),
\end{equation}
where $s_0 = {(c'+d')}/{2}$ and
$$ h_\nu(t, s_0)=
(-1)^{\nu}\frac{1}{\nu!} \prod\limits_{l=1}^\nu  (1-\alpha-l+1)(t-s_0)^{1-\alpha-\nu}, \; g_\nu(s,s_0)=(s-s_0)^{\nu}, $$
with the error estimate
\begin{equation}\label{eqn:errkernel}
    |k(t,s)-\tilde{k}(t,s)|\leq \frac{1}{(a'-d')^{1-\alpha}}  \left(1+\frac{\lambda}{2} \right)   \left(\frac{\lambda}{2+{\lambda}} \right)^{R},
\end{equation}
where $\lambda$ is given in \eqref{admis}.
\end{lem}
\begin{proof}
Since  $ \tau\times \sigma$ are admissible,
we can take the Taylor expansion of the kernel $k(t,s)$ at $s_0$
\begin{align*}
k(t,s) &= \sum_{\nu=0}^\infty  \frac{1}{\nu!} \partial_s^\nu [(t-s_0)^{1-\alpha}] (s-s_0)^{\nu}\\
&= \sum_{\nu=0}^{\infty}  (-1)^{\nu} \frac{1}{\nu!}\prod_{l=1}^\nu (1-\alpha-l+1) (t-s_0)^{1-\alpha-\nu} (s-s_0)^{\nu}.
\end{align*}
By truncating to the term $\nu = R-1$, we obtain the equation \eqref{eqn:Lexpan}.
Denoting  $r= \frac{d'-s_0}{(d'-s_0)+(a'-d')}$,  the  reminder $k(t,s)-\tilde{k}(t,s) $ is estimated by
\begin{equation}\label{eqx001}
    \begin{aligned}
      |k(t,s)-\tilde{k}(t,s)|
 &=
 \left|\sum_{\nu=R}^\infty \frac{1}{\nu!} \partial_s^\nu [(t-s_0)^{1-\alpha}] (s-s_0)^{\nu}\right|
 \leq
\sum_{\nu=R}^\infty \frac{(s-s_0)^{\nu}}{(t-s_0)^{\nu+\alpha-1}}.
    \end{aligned}
\end{equation}
By direct calculation, we have
\be
\frac{s-s_0}{(t-s_0)} \leq \frac{d'-s_0}{(d'-s_0)+(a'-d')}=r,  \quad    \frac{1}{(t-s_0)^{\alpha-1}} \leq \frac{1}{(a'-d')^{\alpha-1}}.
\ee
Consequently, we have
\beq\label{eqx003}
\sum_{\nu=R}^\infty \frac{(s-s_0)^{\nu}}{(t-s_0)^{\nu+\alpha-1}}
=
\sum_{\nu=R}^\infty \frac{1}{(t-s_0)^{\alpha-1}}  \left(\frac{s-s_0}{t-s_0}\right)^\nu
\leq
\sum_{\nu=R}^\infty \frac{r^{\nu}}{(a'-d')^{\alpha-1}}  .
\eeq
By virtue of the admissibility condition \eqref{admis}, the parameter $r$ is  estimated by
\beq\label{eqx006}
r= \frac{ Diam\{\sigma\}/2}{ Diam\{\sigma\}/2 + Dist\{\tau, \sigma\}} \leq  \frac{\lambda}{2+\lambda}.
\eeq
Thus, the summation of the series in \eqref{eqx003} is estimated by
\beq\label{eqx002}
\sum_{\nu=R}^\infty \frac{r^{\nu}}{(a'-d')^{\alpha-1}}   = \frac{r^R}{(a'-d')^{\alpha-1}} \frac{1}{1-r}
\leq \frac{1}{(a'-d')^{\alpha-1}} \frac{2+\lambda}{2} \left(\frac{\lambda}{2+\lambda}\right)^R.
\eeq
Then the conclusion holds by combining \eqref{eqx001}, \eqref{eqx003}  and \eqref{eqx002} .
\end{proof}

\begin{remark}
Instead of expanding the kernel $k(t,s)$ with  respect to $s$, we can also expand the kernel $k(t,s)$ with respect to $t$ at a point $t_0$.
If using an expansion with respect to $t$  at a point $t_0$, then the corresponding remainder is
\begin{equation}\label{eqx004}
 |k(t,s)-\tilde{k}(t,s)|
 =
 \left|\sum_{\nu=R}^\infty \frac{1}{\nu!} \partial_t^\nu [(t_0-s)^{1-\alpha}] (t-t_0)^{\nu}\right|
 \leq
\sum_{\nu=R}^\infty \frac{(t-t_0)^{\nu}}{(t_0-s)^{\nu+\alpha-1}}.
 \end{equation}
Estimations  \eqref{eqx001} and \eqref{eqx004} show how to choose the better expansion:
If $\sup\{|s-s_0| : s\in \sigma \}\leq  \sup\{|t-t_0| : t \in \tau \}$,
the expansion with respect to $s$ is
more favorable, otherwise the expansion with respect to $t$ is better.

Similarly, we can deal with the kernel in the case $t<s$, in this case we have $k(t,s) = (s-t)^{1-\alpha}$ corresponding to the kernel of the right fractional derivative.
\end{remark}

For the sake of simplicity, we denote $h_\nu(t):= h_\nu(t,s_0)$ and $g_\nu(s):= g_\nu(s,s_0)$ if no confusion arises.

\subsection{Low rank approximation of submatrices}
\begin{figure}[!t]
\centering
\includegraphics[width=0.33\textwidth,height=0.33\textwidth]{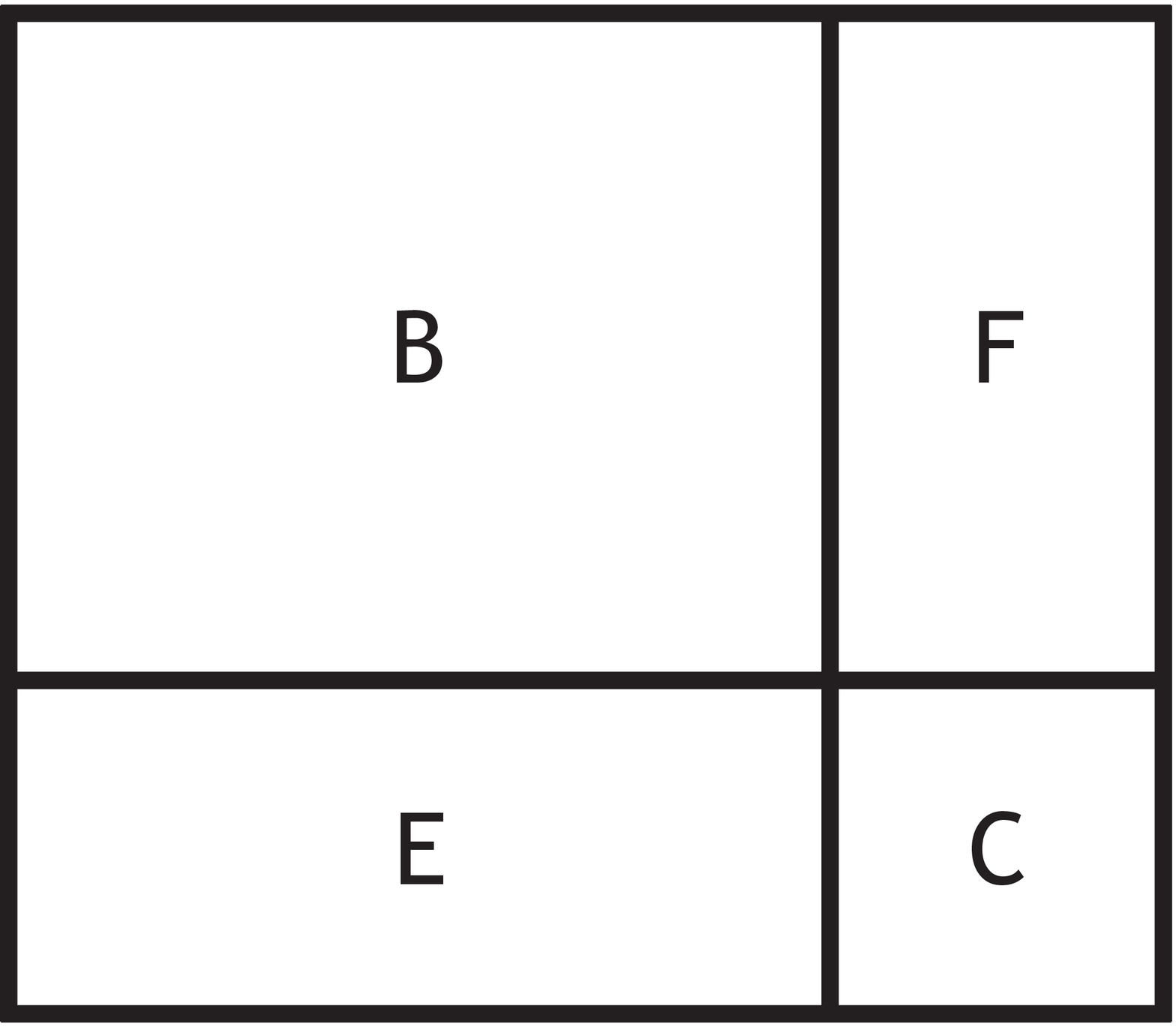}
\quad \quad \quad \quad \quad
\includegraphics[width=0.33\textwidth,height=0.33\textwidth]{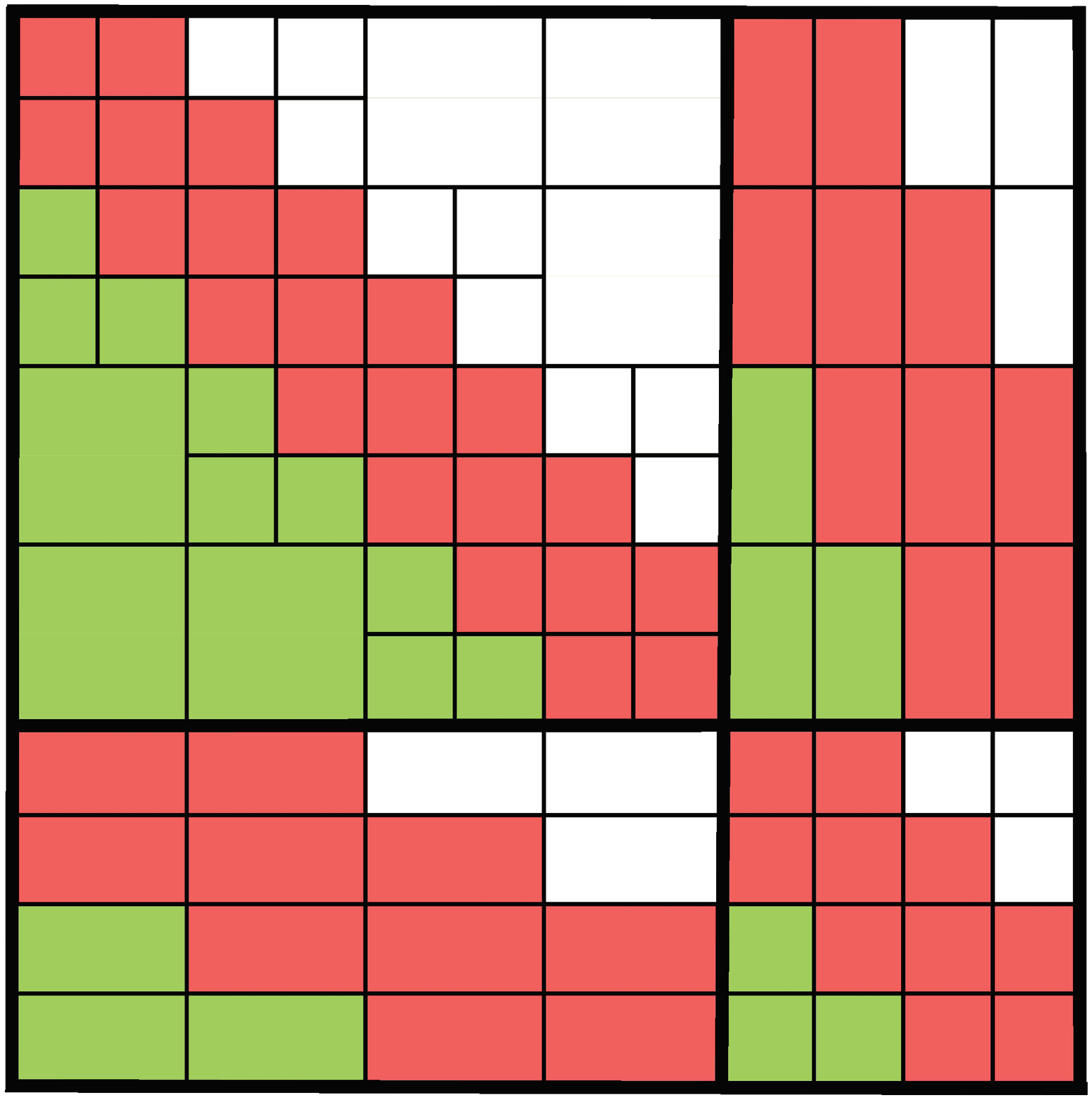}
\caption{\scriptsize Left: Partition of stiffness matrix $S_l$. Right: Green blocks denote low-rank approximation stored in a data sparse format, red blocks stored in a full-matrix representation, white blocks denote zero.
}
\label{fig200} 
\end{figure}

By using the degenerate kernel $\tilde{k}(t, s)$, we can approximate each  submatrix of the matrix $S_l$, namely, $B,\,C,\, E,\, F$, by low rank approximation.
Figure \ref{fig200} shows the schematic diagram of $S_l$ and its H-matrix approximation $H$.

Similar to the admissibility of two domains used for characterizing the local property of the degenerate approximation,
correspondingly, we need to define the admissibility of the index set.
Let $\vartheta:= \{1, \cdots , N\}$ be the index set
and $\vartheta_1,\; \vartheta_2$ be two generic index subsets of $\vartheta$.
Let $S_1$ and $S_2$ be the two corresponding groups of basis function (model basis functions $\{\phi_{i}^k\}_{k=1}^{P_i-1},i=1,\ldots,N$ or nodal basis functions $\{\varphi_{j}\}_{j=1}^{N-1}$) and
\begin{equation*}
    \tau_{\vartheta_1}=\bigcup\limits_{i\in \vartheta_1} \text{supp} \{\phi_i| \phi_i\in S_1\}
, \quad
\tau_{\vartheta_2}=\bigcup\limits_{i\in \vartheta_2} \text{supp} \{\phi_i | \phi_i\in S_2\}
\end{equation*}
be the union of the supports of $S_1$ and $S_2$, respectively.
The sets $ \vartheta_1 \times \vartheta_2$  are said to be admissible if $\tau_{\vartheta_1}\times \tau_{\vartheta_2}$ satisfy the admissibility  condition \eqref{admis}.

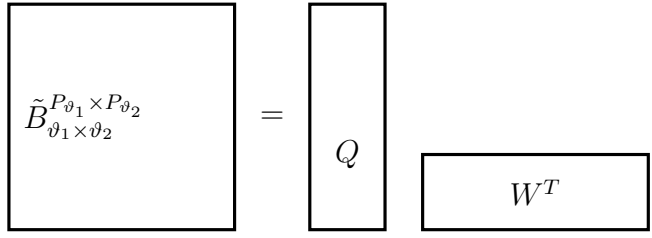
\begin{figure}[!t]
\centering
\begin{tikzpicture}
\draw[black, very thick] (2,0)      rectangle (5, 3);
\draw[black, very thick] (6,0)      rectangle (7, 3);
\draw[black, very thick] (7.5,0)    rectangle (10.5, 1);
\node at (3,1.5) {\large$\tilde{B}_{\vartheta_1 \times \vartheta_2}^{P_{\vartheta_1} \times P_{\vartheta_2}}$};
\node at (5.5,1.5) {\large$=$};
\node  at (6.5,1) {\large$Q$};
\node at (9, 0.5) {\large$W^T$};
\end{tikzpicture}
\caption{\scriptsize Representation of each green block of block B in Figure \ref{fig200} (Right) using a low rank approximation by the block
$\tilde{B}_{\vartheta_1 \times \vartheta_2}^{P_{\vartheta_1} \times P_{\vartheta_2}}$ in terms of matrices $Q$ and $W$.}
\label{fig100}
\end{figure}

If  $\vartheta_1 \times \vartheta_2$ are admissible with respect to two groups of model basis functions  $\{\phi_{i}^k\}_{k=1}^{P_i-1}$ and $\{\phi_{j}^k\}_{k=1}^{P_j-1}$, for $(i,j)  \in \vartheta_1 \times \vartheta_2$, then we can approximate the matrix entries  $B^{pq}_{ij}$  of $B$ by $\tilde{B}^{pq}_{ij}$: Using Lemma 4 of \cite{MaoShen17Adv}, we have
\begin{align*}
  \tilde{B}_{ij}^{pq}
&=
 \frac{1}{\Gamma(2-\alpha)}\int_{x_{i-1}}^{x_i} (\phi_{i}^p(t))'\d t  \int_{x_{j-1}}^{x_j}
\sum_{\nu=0}^{R-1} h_\nu(t) g_\nu(s) ( \phi^{q}_j(s))' \d s \\
&=
\frac{1}{\Gamma(2-\alpha)}\sum_{\nu=0}^{R-1} \int_{x_{i-1}}^{x_i} h_\nu(t) (\phi^{p}_i(t))' \d t
\int_{x_{j-1}}^{x_j} g_\nu(s)  (\phi^{q}_j(s))' \d s   \\
&= \sum_{\nu=0}^{R-1}  Q^p_{i\nu} W_{j\nu}^q,\quad  1\le p\le P_i -1,\, 1\le q\le P_j -1.
\end{align*}
%
%
This means that  the double integral is separated into a multiplication of two single integrals.
More precisely, the submatrix $\tilde{B}_{\vartheta_1 \times \vartheta_2}^{P_{\vartheta_1} \times P_{\vartheta_2}}$
is factorized into
\begin{equation*}
    \tilde{B}_{\vartheta_1 \times \vartheta_2}^{P_{\vartheta_1} \times P_{\vartheta_2}}
=QW^T, \quad  Q\in {\mathbb{R}}^{|P_{\vartheta_1}| \times R},
W\in {\mathbb{R}}^{|P_{\vartheta_2}| \times R},
\end{equation*}
where $ P_{\vartheta_i}=\{ P_m -1 | m\in {\vartheta_i}\}$,
$|P_{\vartheta_i}|=\sum\limits_{ m\in \vartheta_i } P_m- |\vartheta_i |$,
$|\vartheta_i |$ is cardinality of $\vartheta_i$,   $ i=1, 2$,
 and
\begin{equation*}
     Q^p_{i\nu}
=
\frac{1}{\Gamma(2-\alpha)} \int_{x_{i-1}}^{x_i} h_\nu(t) (\phi^{p}_i(t))' \d t,
 \quad
W_{j\nu}^q
= \int_{x_{j-1}}^{x_j}  g_\nu(s)  (\phi^{q}_j(s))' \d s.
\end{equation*}
We show this representation in Figure \ref{fig100}.

Similarly, we have
\begin{equation*}
   \tilde{C}_{\vartheta_1 \times \vartheta_2} =\bar{Q}\,\bar{W}^T, \;
    \tilde{E}_{\vartheta_1 \times \vartheta_2}^{P_{\vartheta_2}}=\bar{Q}\,{W}^T,\;
    \tilde{F}_{\vartheta_1 \times \vartheta_2}^{P_{\vartheta_1}}={Q}\,\bar{W}^T,
\end{equation*}
%
where $\bar{Q}\in {\mathbb{R}}^{|\vartheta_1| \times R},\; \bar{W}\in {\mathbb{R}}^{|\vartheta_2| \times R}$,
\begin{equation*}
    \bar{Q}_{i\nu}
=
\frac{1}{\Gamma(2-\alpha)} \int_{x_{i-1}}^{x_{i+1}} h_\nu(t) (\varphi_i(t))' \d t,
\quad
      \bar{W}_{j\nu} = \int_{x_{j-1}}^{x_{j+1}}  g_\nu(s)  (\varphi_j(s))' \d s.
\end{equation*}
Finally, the  stiffness
matrix $S_l$  can be  approximated by the H-matrix denoted by $H$, i. e.,
\beq\label{approxmat}
S_l\approx H := \left[
\ba{cc}
\tilde{B} & \tilde{F}\\
\tilde{E} & \tilde{C}
\ea
\right].
\eeq
We point out here that for all the elements without using low rank representation we use full-matrix representation, namely, we use the same value as in the original matrix $S$. The structure of matrix $H$ is shown in the Fig. \ref{fig200} (right).

The following theorem provides the elemental error of the approximation  \eqref{approxmat}.
\begin{thm}\label{thm:errest}
Assume $\vartheta_1 \times \vartheta_2$ are admissible and $\tau_{\vartheta_1}=[a',b']$,  $\tau_{\vartheta_2}=[c',d']$,  $d'<a'$,
let $\lambda$ be given as in \eqref{admis}. For $(i,j)\in \vartheta_1 \times \vartheta_2$, we have the approximation errors for $\tilde{B}_{ij}^{pq}$,
$\tilde{E}_{ij}^{q}$, $\tilde{F}_{ij}^{p}$, $\tilde{C}_{ij}$ as follows:
\begin{equation}\label{eq408}
\begin{aligned}
  &|B_{ij}^{pq}-\tilde{B}_{ij}^{pq}| \leq \frac{8h_ic_0}{\Gamma(2-\alpha)}\delta_{\lambda}^{R},
  \quad
  |C_{ij}-\tilde{C}_{ij}| \leq \frac{(h_i+h_{i+1})h_ic_0}{\Gamma(2-\alpha)}\delta_{\lambda}^{R},
  \\
  &|E_{ij}^{q}-\tilde{E}_{ij}^{q}|\leq\frac{2(h_i+h_{i+1})c_0}{\Gamma(2-\alpha)}\delta_{\lambda}^{R},
  \quad
  |F_{ij}^{p}-\tilde{F}_{ij}^{p}| \leq \frac{4h_ic_0}{\Gamma(2-\alpha)}\delta_{\lambda}^{R},
\end{aligned}
\end{equation}
where $\delta_{\lambda} =\lambda/(2+\lambda)$,
$c_0={1}/(|a'-d'|+|d'-s_0|)^{\alpha} \left( 1+\lambda/2 \right)\left( R+1 + \lambda /2 \right)$.
\end{thm}
\begin{proof}
We begin by  estimating  elements of the low rank matrices, i.e., elements of $Q,\, W,\, \bar{Q},\, \bar{W}$.
We first estimate elements of $Q$. Since $h_\nu(t)$ is positive,  by  integrating by parts and applying the mean value theorem for integral, we have
\begin{align*}
   Q^p_{i\nu}
&= \frac{-1}{\Gamma(2-\alpha)} \int_{x_{i-1}}^{x_i} h_\nu'(t) \phi^{p}_i(t) dt
= \frac{-\phi^{p}_i(\xi_i) }{\Gamma(2-\alpha)} \int_{x_{i-1}}^{x_i} h_\nu'(t)  dt \\
&=\frac{-\phi^{p}_i(\xi_i)}{\Gamma(2-\alpha)}  \frac{1}{\nu!} \prod\limits_{l=1}^\nu  (\alpha+l-2)
((x_i-s_0)^{1-\alpha-\nu}- (x_{i-1}-s_0)^{1-\alpha-\nu}),
\end{align*}
where $\xi_i\in [x_{i-1}, x_i]$.
Furthermore, by the mean value theorem, we have
\begin{equation*}
    (x_i-s_0)^{1-\alpha-\nu}- (x_{i-1}-s_0)^{1-\alpha-\nu}
 =
 (1-\alpha-\nu) h_i (\eta_{i}-s_0)^{-\alpha-\nu},
\end{equation*}
where $\eta_i \in [x_{i-1}, x_i]$.
Taking into account  $|\phi^{p}_i(\xi_i)|\leq 2$ noting that $x_i \in [a',b']$, the above two equations yields
\beq\label{eq404}
| Q^p_{i\nu}|
\leq\frac{2|1-\alpha-\nu|h_i}{\Gamma(2-\alpha)}
\left(\frac{1}{|a'-d'|+|d'-s_0|}\right)^{\alpha+\nu}.
\eeq
Then, for the element of $W$, i.e., $W_{j\nu}^q $,
we claim that the following estimate holds:
\beq\label{eq405}
|W_{j\nu}^q|
\leq
4 (d'-s_0)^{\nu}.
\eeq
If $x_j < s_0$, or  $x_{j-1} > s_0$, we can prove \eqref{eq405} by arguing as follows:
by using the similar argument as that for $Q^p_{i\nu}$, we deduce by the mean value theorem that
\begin{equation*}
    W_{j\nu}^q
= -\int_{x_{j-1}}^{x_j}  g_\nu'(s) \phi^{q}_j(s) \d s\\
=
-\phi^{q}_j(\xi_j)   ((x_j-s_0)^{\nu}-(x_{j-1}-s_0)^{\nu} ),
\end{equation*}
where $\xi_j\in [x_{j-1}, x_j]$.
Then, we can obtain \eqref{eq405} since $x_{j-1},x_j \in [c',d']$.
%
If $x_{j-1} \leq s_0 \leq  x_{j} $,  similarly, we have
\begin{align*}
  W_{j\nu}^q
&= -\int_{x_{j-1}}^{x_j}  g'_\nu(s,s_0)  \phi^{q}_j(s) \d s\\
&=
-\nu \left(\int_{x_{j-1}}^{s_0}  (s-s_0)^{\nu-1}  \phi^{q}_j(s) \d s + \int_{s_0}^{x_j}  (s-s_0)^{\nu-1}  \phi^{q}_j(s) \d s\right)\\
&=
- \phi^{q}_j(\xi_j^1)(x_j-s_0)^{\nu}+\phi^{q}_j(\xi_j^2)(x_{j-1}-s_0)^{\nu},
\end{align*}
%
where $\xi_j^1,\xi_j^2\in [x_{j-1}, x_j]$.
Then, we again obtain the estimate \eqref{eq405}.
Next, we estimate $\bar{Q}_{i\nu}$.  Direct computation gives
\begin{equation}\label{eq412}
\begin{aligned}
   \bar{Q}_{i\nu}
=&
\frac{1}{\Gamma(2-\alpha)}\Big(\frac{1}{h_i} \int_{x_{i-1}}^{x_{i}} h_\nu(t)  \d t
-\frac{1}{h_{i+1}} \int_{x_{i}}^{x_{i+1}} h_\nu(t)  \d t \Big)\\
=&
  \frac{1}{\Gamma(2-\alpha)\nu!} \Big(\prod_{l=1}^\nu  (\alpha+l-2)
 \frac{(x_i-s_0)^{2-\alpha-\nu}-(x_{i-1}-s_0)^{2-\alpha-\nu}}{h_i (2-\alpha-\nu)}\\
& - \prod_{l=1}^\nu  (\alpha+l-2)
 \frac{(x_{i+1}-s_0)^{2-\alpha-\nu}-(x_{i}-s_0)^{2-\alpha-\nu}}{h_{i+1}(2-\alpha-\nu)}\Big).
\end{aligned}
\end{equation}
By applying the Taylor expansion, we obtain
\begin{align*}
  (x_{i-1}-s_0)^{2-\alpha-\nu}
=&
(x_i-s_0)^{2-\alpha-\nu} -(2-\alpha-\nu) h_i (x_i-s_0)^{1-\alpha-\nu} \\
&+ {(2-\alpha-\nu)(1-\alpha-\nu)}/{2}\cdot h_i^2 (\eta_i-s_0)^{-\alpha-\nu}, \\
(x_{i+1}-s_0)^{2-\alpha-\nu}
=&
(x_{i}-s_0)^{2-\alpha-\nu} +(2-\alpha-\nu) h_{i+1} (x_{i}-s_0)^{1-\alpha-\nu}\\
&+ {(2-\alpha-\nu)(1-\alpha-\nu)}/{2}\cdot h_{i+1}^2 (\eta_{i+1}-s_0)^{-\alpha-\nu},
\end{align*}
where $\eta_i \in [x_{i-1},x_i]$ and $\eta_{i+1} \in [x_{i},x_{i+1}]$.
Substituting the above two equations into \eqref{eq412} and noting that $x_{i-1},x_i\in [a',b']$, we obtain
\begin{equation}\label{eq411}
\begin{aligned}
  |\bar{Q}_{i\nu}|
\leq &
  \frac{1}{\Gamma(2-\alpha)}
  \frac{|1-\alpha-\nu|}{2}(h_i (\eta_i-s_0)^{-\alpha-\nu}+h_{i+1} (\eta_{i+1}-s_0)^{-\alpha-\nu}) \\
\leq &
  \frac{1}{\Gamma(2-\alpha)}
  \frac{|1-\alpha-\nu|}{2}(h_i+h_{i+1})
  \left(\frac{1}{|a'-d'|+|d'-s_0|}\right)^{\alpha+\nu}.
\end{aligned}
\end{equation}
$\bar{W}$ can be estimated in a similar way, for instance,
\begin{equation}\label{eq415}
\begin{aligned}
|\bar{W}_{j\nu}| = &
\left| \frac{1}{h_j}\int_{x_{j-1}}^{x_{j}}  g_\nu(s)   \d s
-
\frac{1}{h_{j+1}}\int_{x_{j}}^{x_{j+1}}  g_\nu(s)  \d s \right| \\
= &
\left| \frac{(x_{j}-s_0)^{\nu+1}-(x_{j-1}-s_0)^{\nu+1}}{h_j (\nu+1)}
-
\frac{(x_{j+1}-s_0)^{\nu+1}-(x_{j}-s_0)^{\nu+1}}{h_{j+1}  (\nu+1) }\right| \\
\leq &
|\eta_{j}-s_0|^{\nu}  + |\eta_{j+1}-s_0|^{\nu}
\leq
2 (d'-s_0)^{\nu},
\end{aligned}
\end{equation}
where $\eta_i \in [x_{i-1},x_i]$ and $\eta_{i+1} \in [x_{i},x_{i+1}]$.

Now we turn to estimate the elements of $B,\,C,\, E,\,F$.
Denote $r= \frac{d'-s_0}{(d'-s_0)+(a'-d')}$.
By the equations \eqref{eq404} and \eqref{eq405}, we deduce
\begin{equation}\label{eqx008}
\begin{aligned}
&\left|B_{ij}^{pq}-\tilde{B}_{ij}^{pq} \right| = \left| \sum_{\nu=R}^\infty  Q^p_{i\nu} W_{j\nu}^q \right| \\
\leq&
\sum_{\nu=R}^\infty \frac{8h_i(\alpha+\nu-1)}{\Gamma(2-\alpha)}
\left(\frac{1}{|a'-d'|+|d'-s_0|}\right)^{\alpha+\nu}
  (d'-s_0)^{\nu}\\
  \leq&
\frac{8h_i}{\Gamma(2-\alpha)} \left(\frac{1}{|a'-d'|+|d'-s_0|}\right)^{\alpha}\sum_{\nu=R}^\infty   (\nu+1)r^{\nu} .
\end{aligned}
\end{equation}
Noting  that $r\leq \frac{\lambda}{2+\lambda} <1$ in  \eqref{eqx006}, the sum of series  appearing in the above inequality is estimated by
 \beq\label{eqx007}
 \sum_{\nu=R}^\infty (\nu+1) r^{\nu} =   \frac{R+1-rR}{(1-r)^2}   r^{R} \leq \left(1+\frac{\lambda}{2}\right) \left(R+1+\frac{\lambda}{2}\right)  \left(\frac{\lambda}{2+\lambda}\right)^{R}.
 \eeq
The first  estimate of  \eqref{eq408} holds by combing  \eqref{eqx008} and \eqref{eqx007}.
Then, we deduce from \eqref{eq411}  and  \eqref{eq415} that
\begin{equation}\label{eqx009}
\begin{aligned}
  &|C_{ij}-\tilde{C}_{ij}|
=| \sum_{\nu=R}^\infty  \bar{Q}_{i\nu} \bar{W}_{j\nu}|\\
\leq&
\sum_{\nu=R}^\infty \frac{|1-\alpha-\nu|(h_i+h_{i+1})}{\Gamma(2-\alpha)}
  \left(\frac{1}{|a'-d'|+|d'-s_0|}\right)^{\alpha+\nu}
  (d'-s_0)^{\nu}
  \\
\leq&
\frac{h_i+h_{i+1}}{\Gamma(2-\alpha)} \left(\frac{1}{|a'-d'|+|d'-s_0|}\right)^{\alpha}
 \sum_{\nu=R}^\infty  (\nu+1)r^{\nu}.
\end{aligned}
\end{equation}
Therefore, the second estimate of  \eqref{eq408} follows from the \eqref{eqx007} and \eqref{eqx009}.
By using the same argument, we can obtain the third and fourth estimates of  \eqref{eq408}.
\end{proof}

\section{Fast solver for the linear system}\label{sec:Fsolver}
Now we turn to solve the linear system \eqref{lsys} by developing a fast solver. Specifically, we first employ the Hierarchical LU (H-LU) decomposition to solve the H-matrix approximation problem, and then solve the  linear system \eqref{lsys} by using the H-matrix approximation as a preconditioner.

\subsection{Solving the H-matrix approximation problem using H-LU decomposition}
\begin{figure}[!t]
\centering
\includegraphics[width=0.92\textwidth]{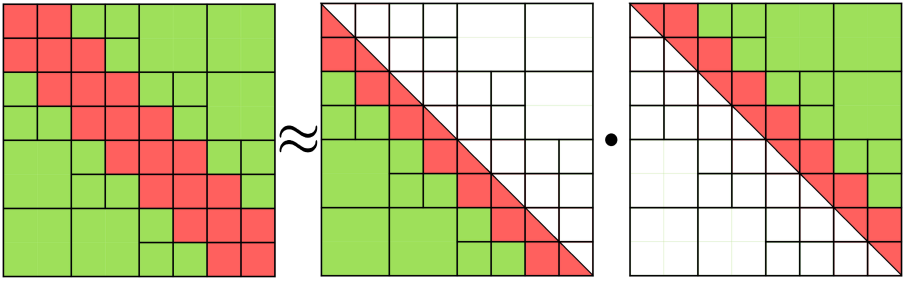}
\caption{\scriptsize H-LU decomposition of $\tilde{A}$. Green blocks denote low-rank approximation stored in a data sparse format,  red blocks stored in a full-matrix representation,  and white blocks denote zero.}
\label{fig:HLU} 
\end{figure}

Instead of solving the original linear system \eqref{lsys}, we solve, in this subsection, the H-matrix approximation problem
\beq\label{approxsys}
\tilde{A}X = G, \text{ where }\tilde{A}:= \rho M + H.
\eeq
In particular, we apply H-LU  decomposition~\cite{borm2005hierarchical} to solve the above system.
There are two main steps to solve the H-matrix approximation problem \eqref{approxsys}. The first step is to implement the H-LU decomposition.
The objective of H-LU decomposition for $\tilde{A}$ is to obtain the decomposition of form
\begin{equation*}
   \tilde{A} \approx L_{H}U_{H}
\end{equation*}
with a prescribed tolerance $Tol_{HLU}$, where $L_{H}$ is a  lower triangular matrix with ones on the diagonal and $U_{H}$ is a  upper triangular matrix.
Moreover, both $L_{H}$ and $U_{H}$ are stored in H-matrix format, see Figure \ref{fig:HLU}.
The second step is to solve a hierarchical lower triangle system $L_{H}Y = G$ using a forward substitution, and then to solve a hierarchical upper triangle system $U_{H}X = Y$ using a backward substitution.
More details can be found in \cite[Section 5.2.3]{borm2005hierarchical}.

The cost of the forward  and the backward substitutions are $O(R N_d  \log N)$.
The cost of H-LU decomposition is $O(R^2 N_d  \log^2N) +O(R^3 N_d  \log N)$, see \cite{borm2005hierarchical}.
Obviously, the cost of H-LU decomposition is dominant, and then the total cost of solving \eqref{approxsys} using H-LU decomposition is $O(R^2 N_d  \log^2N) +O(R^3 N_d  \log N)$.

\subsection{Using the H-matrix approximation system as a preconditioner}
By using H-LU decomposition, we reduce the computation cost to $O(R^2 N_d \log^2N)+O( R^3 N_d \log N)$ for solving the system \eqref{approxsys}. However, it is at the cost of losing accuracy 
and the system $\tilde{A}$ is ill-conditioned. 
Next, we solve the original system \eqref{lsys} by using the H-matrix approximation system \eqref{approxsys} as a preconditioner.
More precisely, we solve the following preconditioned system
\beq\label{precondsys}
\tilde{A}^{-1}AX=\tilde{A}^{-1}G
\eeq
using the BICGSTAB method.

To solve the preconditioned system \eqref{precondsys} with the BICGSTAB method, we need the evaluation of the matrix-vector multiplication $AX$. For a non-uniform mesh, this requires the computational cost to be $O(N_d^2)$, then the total computational cost is $O(N_d^2)+ O(R^2 N_d \log^2  N)+O(R^3 N_d  \log N)$. However, let the polynomial degree in each element be a constant, i.e., $P_i=P$, then for a uniform mesh or a geometric mesh, we can further reduce the computational cost to $O(P^2 N\log N)+ O(R^2 N_d \log^2 N )+O(N_d R^3 \log N) $, where $N$ is the number of elements.

\subsection{Fast evaluation of $AX$ on a uniform or a geometric mesh}\label{sec:toeplitz}

We introduce in this subsection the fast evaluation of matrix-vector multiplication $AX$, where $A = \rho M +S$, when using a uniform mesh or a geometric mesh, i.e., $\hat{q}=\frac{h_{i+1}}{h_i}, i=1, 2, \cdots, N-1 $ is a constant.
The evaluation of the matrix-vector multiplication $MX$ is negligible compared with the one of $SX$ since $M$ is a sparse matrix. Moreover, since $S = \theta S_l +(1-\theta) S_l^T$, hence, we only show the evaluation of $S_l X$. The evaluation of $S_l^T X$ can be implemented in a similar way. Assume that $X$ is divided by $X = [x_1,x_2]^T$ corresponding to the partition of $S_l$. Then we have
\begin{equation*}
    S_l X= \left[
\ba{cc}
{B} & {F}\\
{E} & {C}
\ea
\right]
\left[
\ba{c}
x_1\\
x_2
\ea
\right]=
\left[\ba{c}
B x_1+F x_2\\
E x_1+C x_2
\ea
\right].
\end{equation*}

Now we show how to compute $Bx_1$ at the cost of $O(P^2 N\log N)$.
We need to rearrange the elements of matrix $B$ as $\hat{B}$,
\begin{equation*}
\hat{B}:=\left[
\begin{array}{ccc}
\boxed{B^{11}} & \cdots   & \boxed{B^{1,P-1}}   \\
\vdots          & \ddots & \vdots     \\
\boxed{B^{P-1,1}} & \cdots & \boxed{ B^{P-1,P-1}}\\
\end{array}
\right],
\end{equation*}
where blocks are given by
${B^{pq}}
=\{B_{ij}^{pq}\}_{i,j=1}^{N}$,  $p, q =1,\cdots, P-1$,   and rearrange the vector $x_1$ as $\hat{x}_1$,
\begin{equation*}
    \hat{x}_1=\big[\hat{u}_{1}^{1}, \cdots, \hat{u}_{N}^{1}; \cdots ; \hat{u}_{1}^{P-1}, \cdots, \hat{u}_{N}^{P-1}\big].
\end{equation*}
Let
\begin{equation}\label{eqn:Dh}
    D_h =
\left[
\frac{h_1}{2}, \frac{h_2}{2}, \cdots,  \frac{h_{N}}{2}
\right]^T, \quad D_{h,\alpha} = D_h^{-\alpha}.
\end{equation}
The block $\hat{B}^{pq}$ has form of
\beq\label{eq420}
\hat{B}^{pq}=\mathrm{diag}(D_{h,\alpha})
\,T_B \,
\mathrm{diag}(D_{h}),
\; \text{where }
T_B=
\left[
\ba{cccc}
a^{pq}_0&&&\\
a^{pq}_1&a^{pq}_0&&\\
\vdots&&\ddots&\\
a^{pq}_{N-1}&\cdots &a^{pq}_1& a^{pq}_0
\ea
\right]
\eeq
is a Toeplitz matrix with
\begin{align*}
a^{pq}_0=& \gamma_0\int_{-1}^1 (\psi_p(x))' \d x \frac{d}{d x} \int_{-1}^x  (x- t)^{1-\alpha}  \psi_q(t) \d t  , \\
a^{pq}_k=&\gamma_0\int_{-1}^1 (\psi_p(x))' \d x \frac{d}{d x} \int_{-1}^1  (x-\hat{q}^k t+ 1+2\hat{q}+\cdots+2\hat{q}^{k-1}+\hat{q}^k)^{1-\alpha}  \psi_q(t)    \d t
\end{align*}
%
for $k=1, \cdots, N-1$, where $\gamma_0={1}/{\Gamma(2-\alpha)}$ and $\psi_j(x),j\ge 0$ are defined in \eqref{RefBasis}.

It is known that the Toeplitz matrix $T_B$ can be embedded into a $2N \times 2N$ circulant matrix $C_{2N}$  \cite{bottcher2012introduction,gray2006toeplitz} as follows
\begin{equation*}
C_{2N}= \left[
\ba{cc}
{T_B} & {T_B'}\\
{T_B'} & {T_B}
\ea
\right],
\text{ where }
T'_B=
\left[
\ba{cccc}
0&a^{pq}_{N-1}&\cdots& a^{pq}_1\\
0&0&a^{pq}_{N-1}&\vdots\\
\vdots&&\ddots&a^{pq}_{N-1}\\
0&\cdots &0& 0
\ea
\right].
\end{equation*}
The circulant matrix $C_{2N}$ can be decomposed as  \cite{davis1979circulant}
\begin{equation*}
C_{2N}=F^{-1}_{2N}\, \mathrm{diag}(F_{2N}V_0)\, F_{2N},
\end{equation*}
where $V_0$ is the first column vector of $C_{2N}$ and $F_{2N}$ is the $2N \times 2N $ discrete Fourier transform matrix
\begin{equation*}
F_{N}(j,l)= \frac{1}{\sqrt{N}} \exp\left( -\frac{2\pi {\rm i}jl}{N}  \right),\quad  0\leq j, l\leq N-1,\;  {\rm i}=\sqrt{-1}.
\end{equation*}
Let $x_{2N}$ be a column vector of size $2N$, then the matrix-vector multiplication $F_{2N}x_{2N} $ can be carried out in $O(2N\log(2N)) = O(N\log N) $ operations via the fast Fourier transform.  Thus, $T_B x_{2N} $ can be carried out in $O(N\log N) $ operations, and  hence, $Bx$ can be evaluated in $O(P^2N\log N)$ operations.

The submatrices $C, E , F$ can be decomposed  similarly as in \eqref{eq420} (see Appendix \ref{sec:apd:mf}), and the corresponding computational costs are
$O(N\log N)$, $O(PN\log N)$ and $O(PN\log N)$, respectively.
Thus, the total computational cost of $S_l X$ is $O(P^2 N\log N)$.

\section{Numerical Tests}\label{sec_num}
In this section, we present several numerical tests to verify the error estimates, compare the accuracy, CPU time and condition number of the original system \eqref{lsys} ($AX=G$), the H-matrix approximation system \eqref{approxsys} ($\tilde{A}X = G$) and the preconditioned system \eqref{precondsys} ( $\tilde{A}^{-1}AX=\tilde{A}^{-1} G$).

We point out here that the original system $AX=G$ is solved by using the BICGSTAB iterative method or LU decomposition method, while the H-matrix approximation system $\tilde{A}X = G$ is solved by using the H-LU decomposition method proposed in the previous section, and the preconditioned system $\tilde{A}^{-1}AX=\tilde{A}^{-1} G$ is solved by using preconditioned BICGSTAB method, where the preconditioner $\tilde{A}X = G$ is again solved by the H-LU decomposition method. Moreover, the tolerance of the BICGSTAB for $AX=G$ or the preconditioned BICGSTAB for $\tilde{A}^{-1}AX=\tilde{A}^{-1} G$ is set to $10^{-13}$; the tolerance of the H-LU decomposition for solving the H-matrix approximation system is $10^{-13}$ while for the preconditioner is $10^{-3}$.

\begin{exam}\label{ex:uex}
Let $\theta = 1$. Suppose the  exact solution for the equation \eqref{Model}  is
    $$u(x) = (x-a) - (b-a)^{1-\gamma}  (x-a)^{\gamma}$$
%
with $a=0,\, b=10,\, \gamma=0.8$. In this case, we have $c_1 = c_2 = 0$.
\end{exam}

\begin{figure}[!t]
\centering
\includegraphics[width=0.54\textwidth,height=0.45\textwidth]{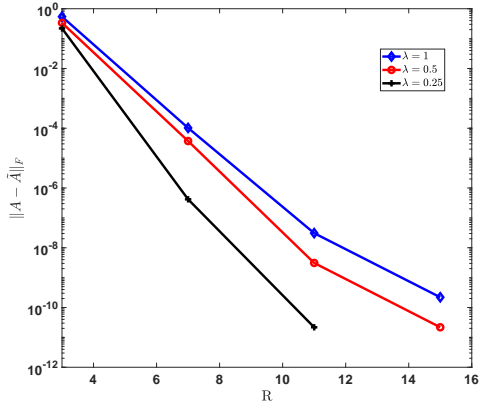}
\caption{\scriptsize Convergence of the error of the matrix $A-\tilde{A}$ in Frobenius norm with respect to the value of rank $R$. $N=600, P=3, \alpha=1.6, \rho = 0, \theta = 1$.}
\label{fig:matrixerror} 
\end{figure}

We begin by considering a graded mesh, namely, $x_i= (\frac{i}{N})^{\tilde{q}}(b-a), i=0,1,\ldots,N$ with $\tilde{q}=5$ and $N$ is the number of elements.
We first verify the error estimate \eqref{eq408} for the matrix entries. To do this, we show the convergence of the errors of $A-\tilde{A}$ in Frobenius norm with respect to the value of rank $R$ for different values of $\lambda$ given in \eqref{admis}. We can see in Figure \ref{fig:matrixerror} that the errors of $\|A-\tilde{A}\|_F$ decay exponentially for all values of $\lambda = 1,0.5,0.25$. This verifies our theoretical analysis and Theorem \ref{thm:errest}. In the following tests, we set $\lambda=1$.

\begin{figure}[!t]
\centering
\includegraphics[width=0.48\textwidth,height=0.4\textwidth]{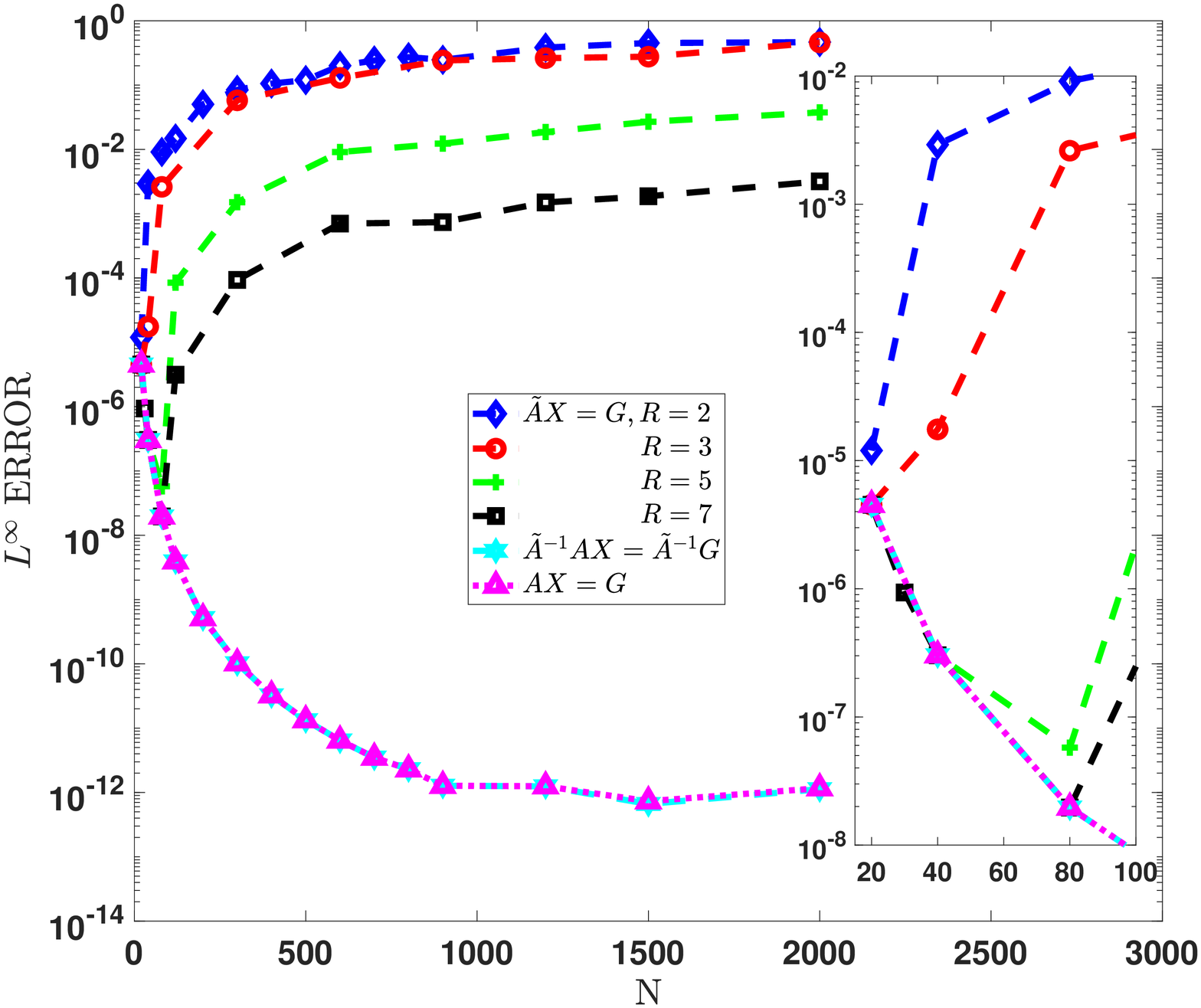}
\includegraphics[width=0.48\textwidth,height=0.4\textwidth]{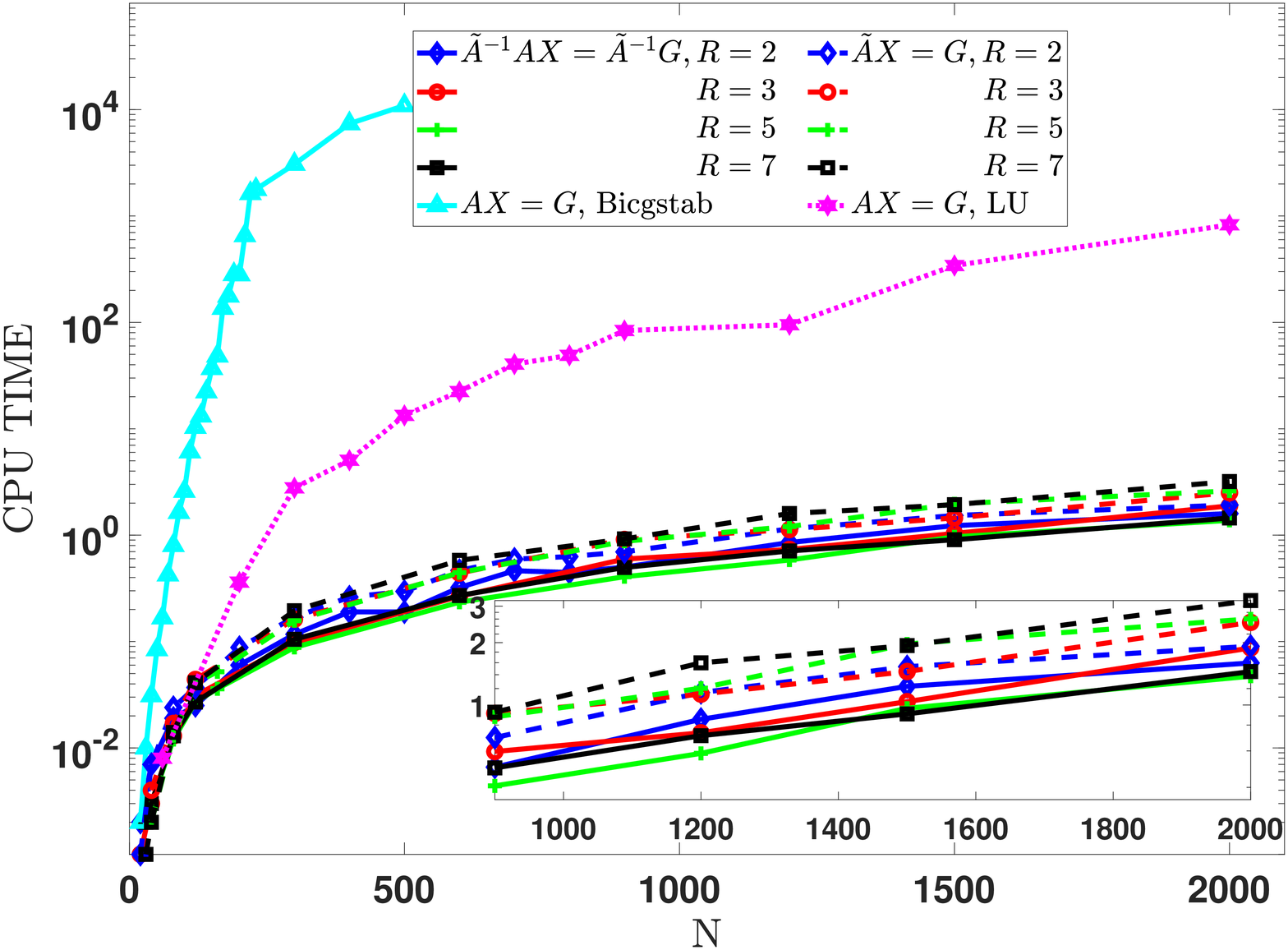}
\caption{\scriptsize Example \ref{ex:uex}:  Accuracy and cost for $h$-refinement. $\rho=0,\,\theta = 1,\, P=3$.  Left: $L^{\infty}$-errors versus number of elements $N$. The preconditioned system has the same accuracy as the original system whereas the H-matrix approximation diverges although initially seems to converge for $R\ge 2$ (see inset).  Right: CPU time versus number of elements $N$. The cases on the left column of the legend are represented by solid lines while the cases on the right column are represented by dash lines. The inset shows a zoom-in plot for clarity.
}
\label{fig3} 
\end{figure}

\begin{figure}[!t]
\centering
\includegraphics[width=0.48\textwidth,height=0.4\textwidth]{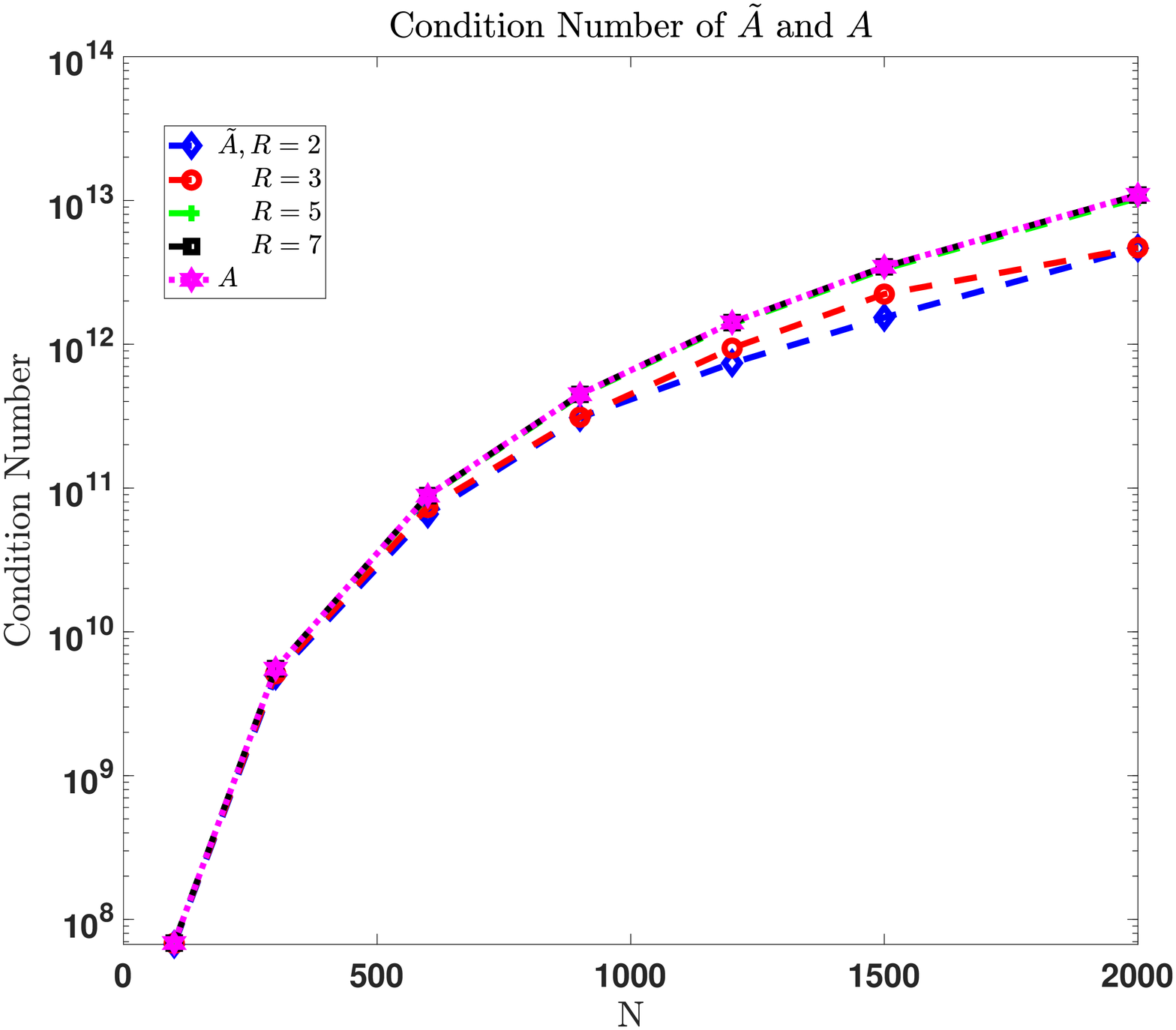}
\includegraphics[width=0.48\textwidth,height=0.4\textwidth]{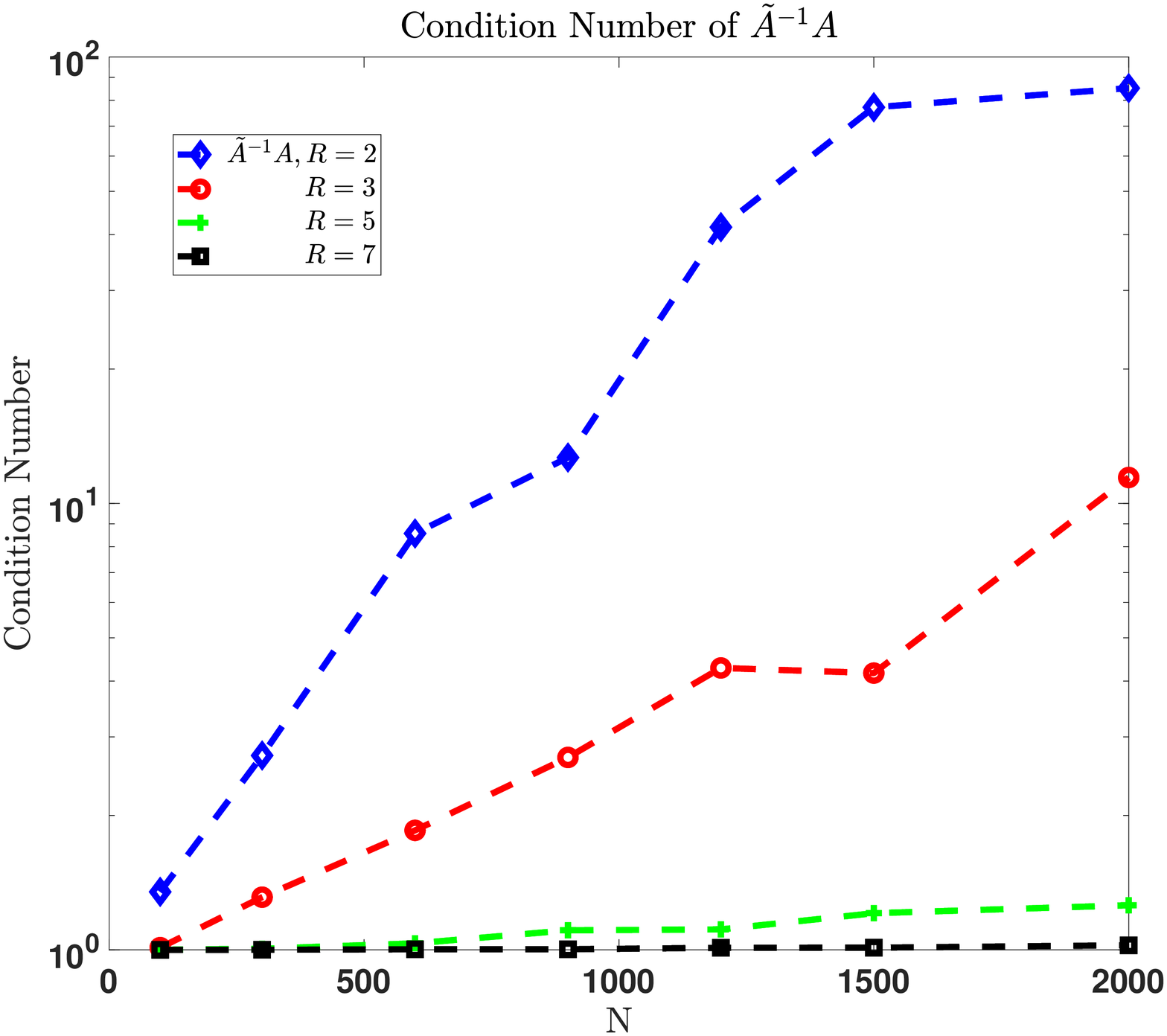}
\caption{\scriptsize Example \ref{ex:uex}:  Condition number for $h$-refinement for different values of rank $R$. $\rho=0,\,\theta = 1,\, P=3$. Left: condition number of $\tilde{A}$ and $A$; right: condition number of $\tilde{A}^{-1}A$.
}
\label{fig32} 
\end{figure}

For the three aforementioned systems,
we now show convergence results for $h$-refinement, namely, the $L^\infty$-error versus the number of element $N$ by fixing the degree of polynomial of each element to be $P_i = P = 3$, in the left plot of Figure \ref{fig3} with $\rho= 0$.
We observe that we can obtain high accuracy in terms of the $L^\infty$-error for the original system. For the H-matrix approximation system, the errors decay similarly as the one of the original system in a certain range of $N$ depending on the value of the rank. However, the errors do not converge and, even worse, they begin to deteriorate as $N$ increases. This is because the truncation error of the H-matrix approximation accumulates as the number of element increases. Then, by using the H-matrix approximation $\tilde{A}$ as a preconditioner, we obtain the same convergence as the original system for the preconditioned system.
We then compare the CPU time (in seconds) for solving the three systems. The results are shown in the right plot of Figure \ref{fig3}. Obviously, the CPU time of solving the original system is much higher than that of solving the H-matrix approximation system. However, after applying the preconditioner, the CPU time is significantly reduced.
Furthermore, to gain some insight, we present the condition numbers of the original system and the approximation system in the left of Figure \ref{fig32} and similarly for the preconditioned system with different values of rank $R$ in the right of Figure \ref{fig32}. We see that the condition numbers of the original system and the approximation system grow very fast. However, the condition number of the preconditioned system is significantly reduced, moreover, the condition number is almost a constant close to 1 when using high value of rank ($R = 7$). We also present the number of iterations for the preconditioned system with different values of rank $R$ in Table \ref{table1}. We can see that the number of iterations does not increase if $R = 7$, which is consistent with our previous observation. This means that the proposed preconditioner is optimal for large rank.

Overall, the preconditioned system possesses both the advantages of the SEM approximation and the H-matrix approximation. In particular, the preconditioned system can be solved efficiently while retaining the high accuracy of the SEM approximation.

\begin{table}[H]
\begin{center}
\caption{Example \ref{ex:uex}:  Number of iterations for solving $\tilde{A}^{-1}AX=\tilde{A}^{-1}G$ using BICGSTAB with respect to $h$-refinement,  $\rho=0,\,\theta = 1,\, P = 3$.}
\begin{tabular}{c|c|c|c|c}
\hline
\hline
   & $R=2$ & $R=3$& $R=5$& $R=7$ \\
 \hline
 $N=100$ & 6 & 6&4&4 \\
  $N=400$& 10 & 8&4&4\\
    $N=800$&12&10&6&4\\
  $N=1500$&16&12&6&4\\
 $N=2000$&16&16&6&4\\
 \hline
\end{tabular}
\end{center}
\label{table1}
\end{table}

\begin{figure}[!t]
\centering
\includegraphics[width=0.48\textwidth]{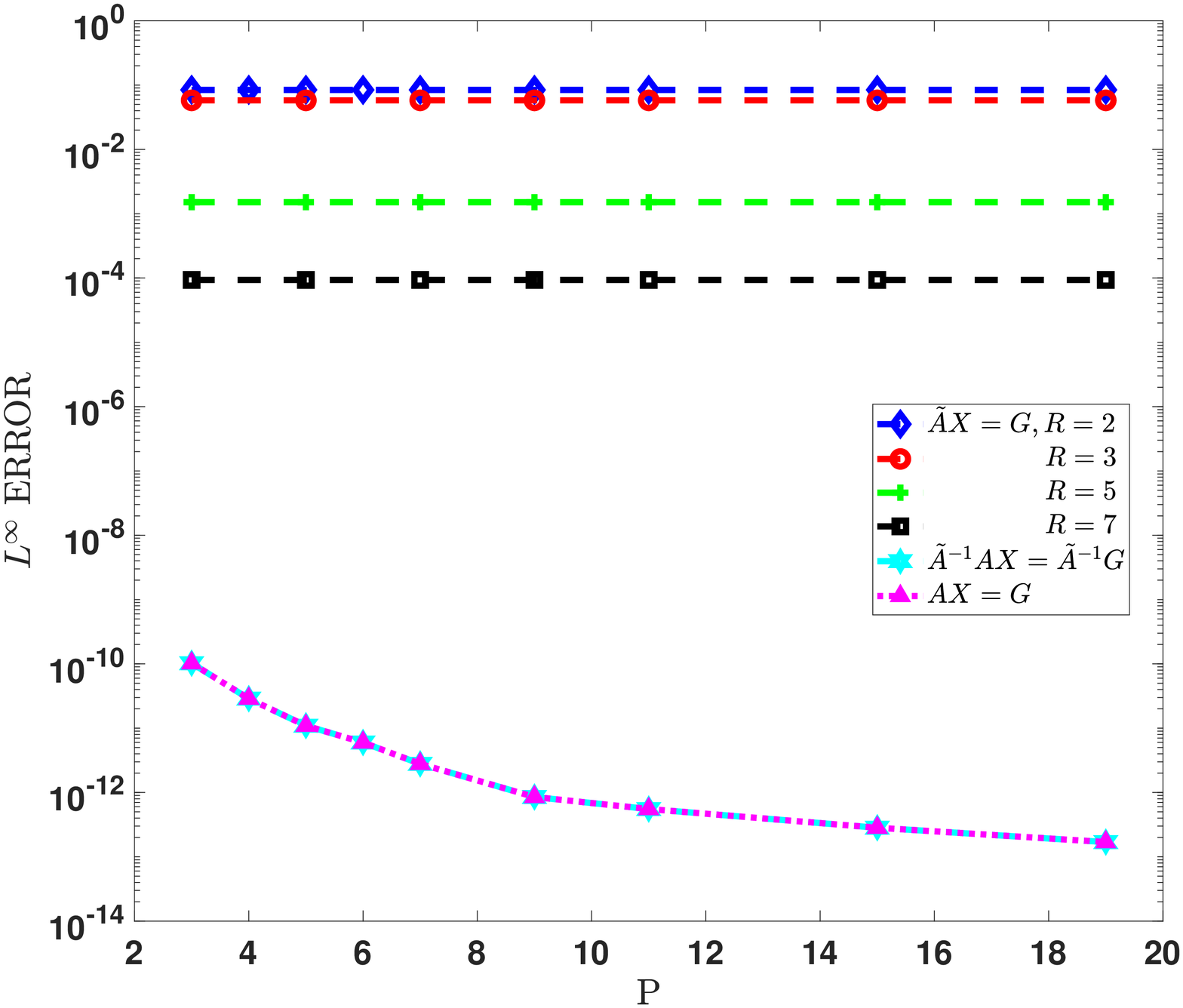}
\includegraphics[width=0.48\textwidth]{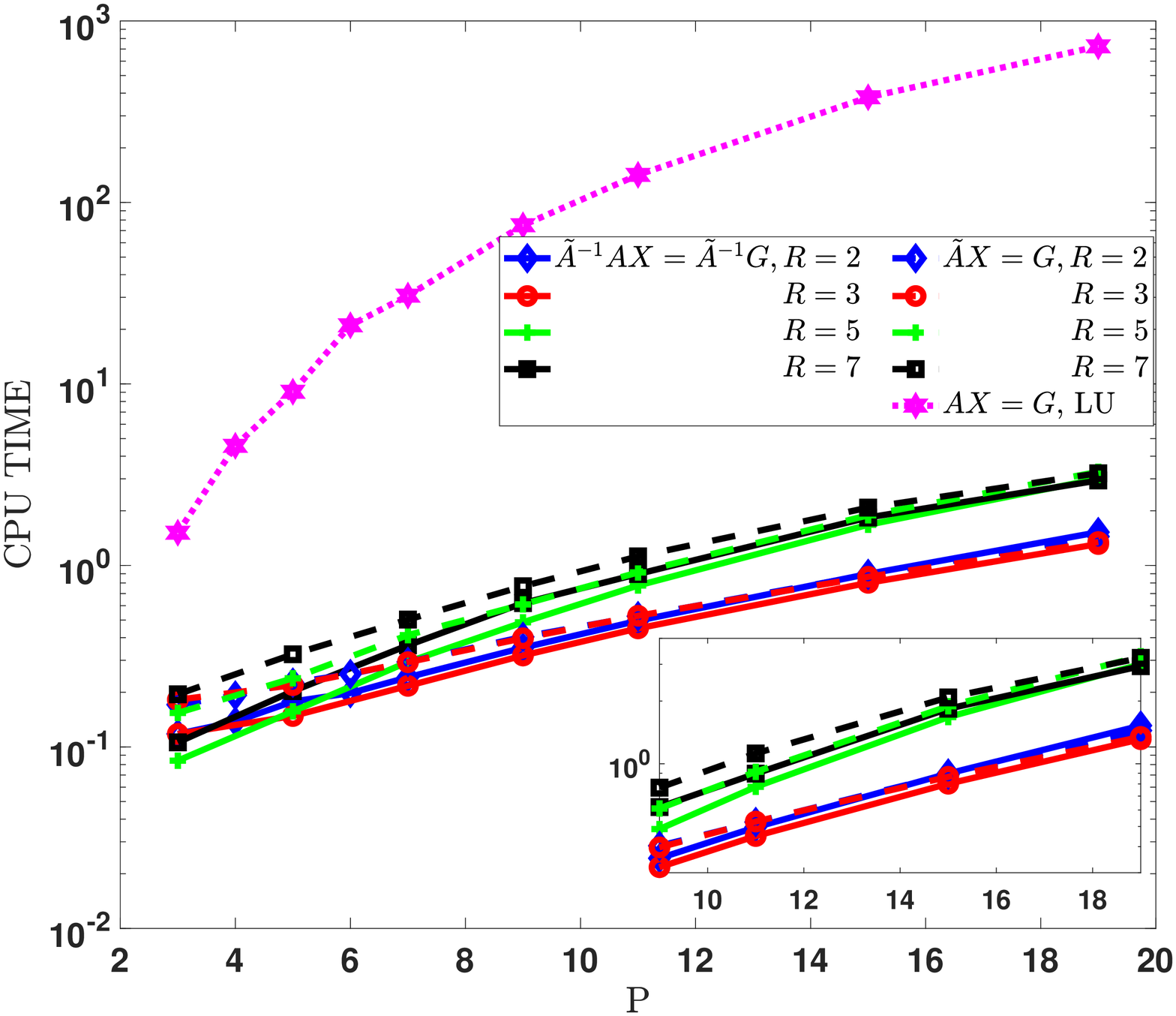}
\caption{\scriptsize Example \ref{ex:uex}:  Accuracy and cost for $p$-refinement. $\rho=0,\,\theta = 1,\, N=300$.  Left: $L^{\infty}$-errors versus degrees of polynomial $P$. The preconditioned system has the same accuracy as the original system whereas the H-matrix approximation does not improve with $p$-refinement. Right: CPU time versus degrees of polynomial $P$. The cases on the left column of the legend are represented by solid lines while the cases on the right column are represented by dash lines. The inset shows a zoom-in plot for clarity.
}
\label{fig4} 
\end{figure}

\begin{figure}[!t]
\centering
\includegraphics[width=0.48\textwidth]{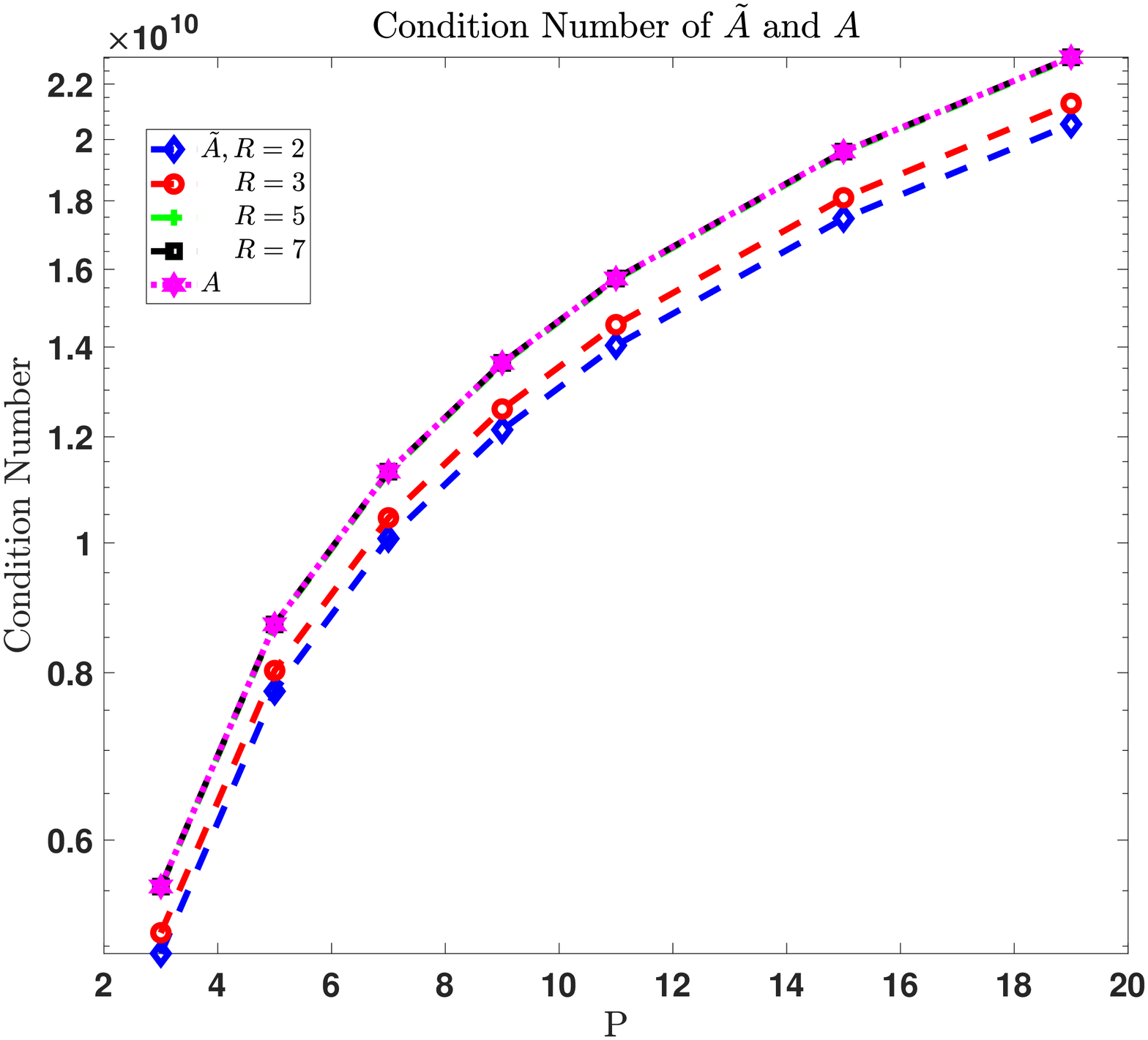}
\includegraphics[width=0.48\textwidth]{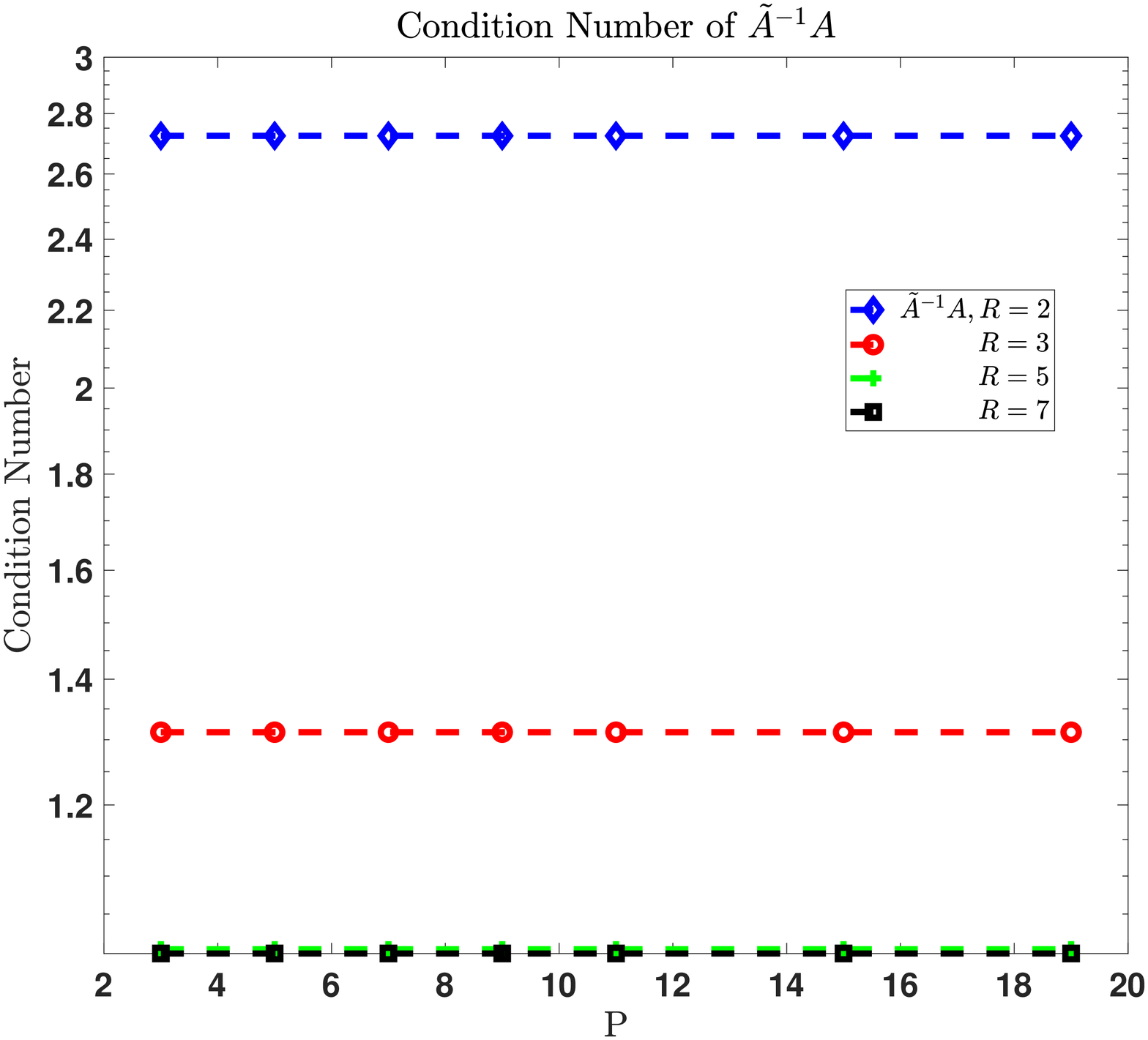}
\caption{\scriptsize Example \ref{ex:uex}:
Condition number for $p$-refinement for different values of rank $R$. $\rho=0,\,\theta = 1,\, N=300$. Left: condition number of $\tilde{A}$ and $A$; right: condition number of $\tilde{A}^{-1}A$.
}
\label{fig42} 
\end{figure}

Then, for the $p$-refinement, we also show the convergence of $L^\infty$-error, CPU time in Figure \ref{fig4}, the condition number in Figure \ref{fig42}, and the number of iterations for the preconditioned system in Table \ref{table2}.
In this case, we set $\rho =0$ and the number of elements to be $N = 300$.
A similar conclusion can be made. Specifically, we can obtain high accuracy for the original system but it is time consuming, while we can efficiently solve the H-matrix approximation system but lose high accuracy. Again, by solving the preconditioned system, we can obtain high accuracy at a much lower computational cost. We point out here that for the $p$-refinement the convergence of the H-matrix approximation is slightly different from the one of $h$-refinement, i.e., when the number of elements is larger than a critical value, the $L^\infty$-error neither decays nor grows for the $p$-refinement while the $L^\infty$-error grows for the $h$-refinement. This means that the $p$-refinement does not suffer from the accumulation of the H-matrix approximation error. The second difference of these two refinements is that the condition number of the preconditioned system is almost constant for any fixed value of rank for the $p$-refinement while the one for $h$-refinement increases monotonically with respect to $N$ with small value of rank.

\begin{table}[H]
\begin{center}
\caption{Example \ref{ex:uex}:  Number of iterations for solving $\tilde{A}^{-1}AX=\tilde{A}^{-1}G$ using BICGSTAB with respect to $p$-refinement,  $\rho=0,\,\theta = 1,\, N = 300$.}
\begin{tabular}{c|c|c|c|c}
\hline
\hline
   & $R=2$ & $R=3$& $R=5$& $R=7$ \\
 \hline
 $P=3$ & 8 & 8&6&4 \\
  $P=5$& 8 & 8&6&4\\
    $P=9$&8&8&6&4\\
  $P=13$&8&8&6&4\\
 $P=19$&8&8&6&4\\
 \hline
\end{tabular}
\end{center}
\label{table2}
\end{table}

\begin{figure}[!t]
\centering
\includegraphics[width=0.48\textwidth]{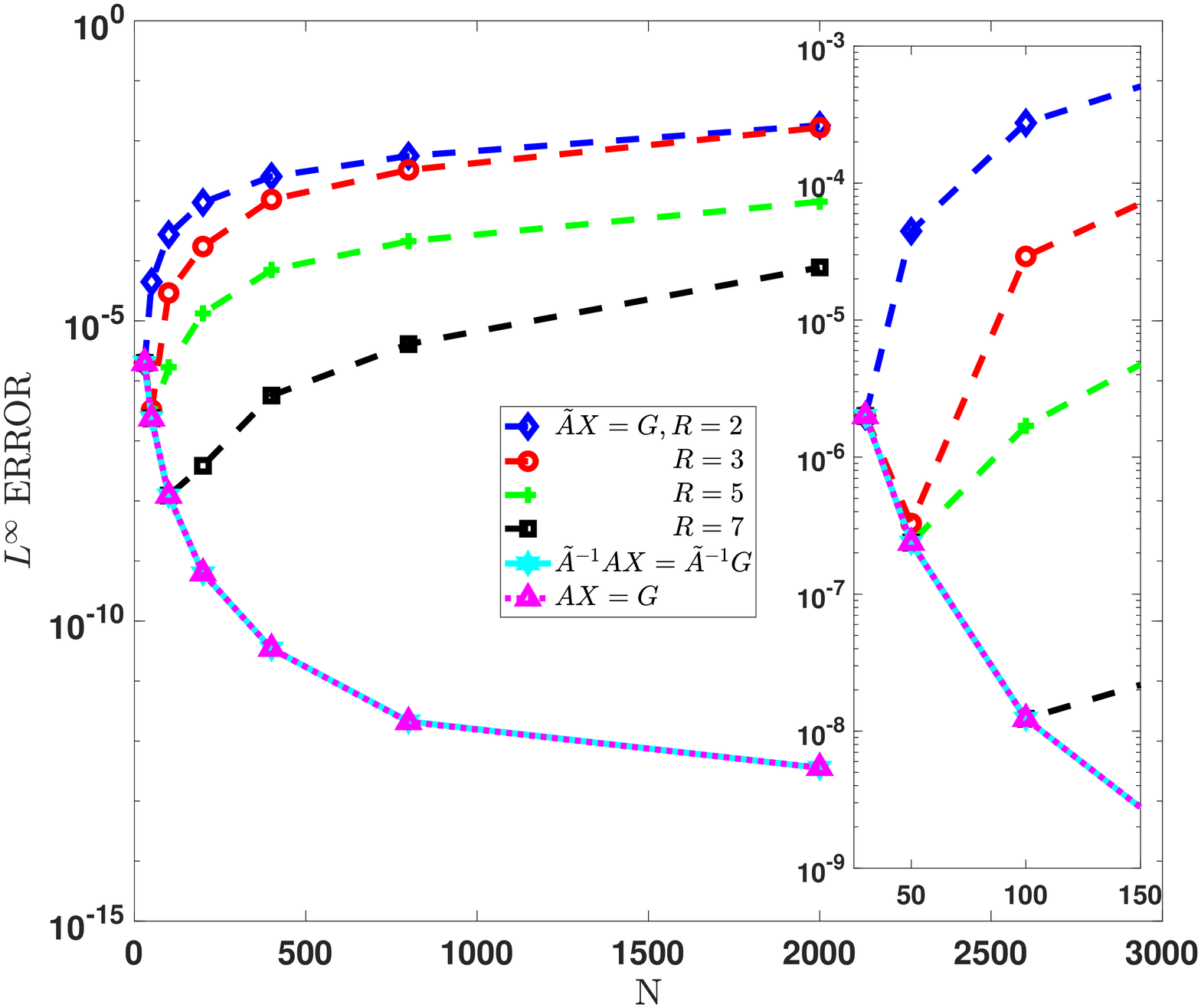}
\includegraphics[width=0.48\textwidth]{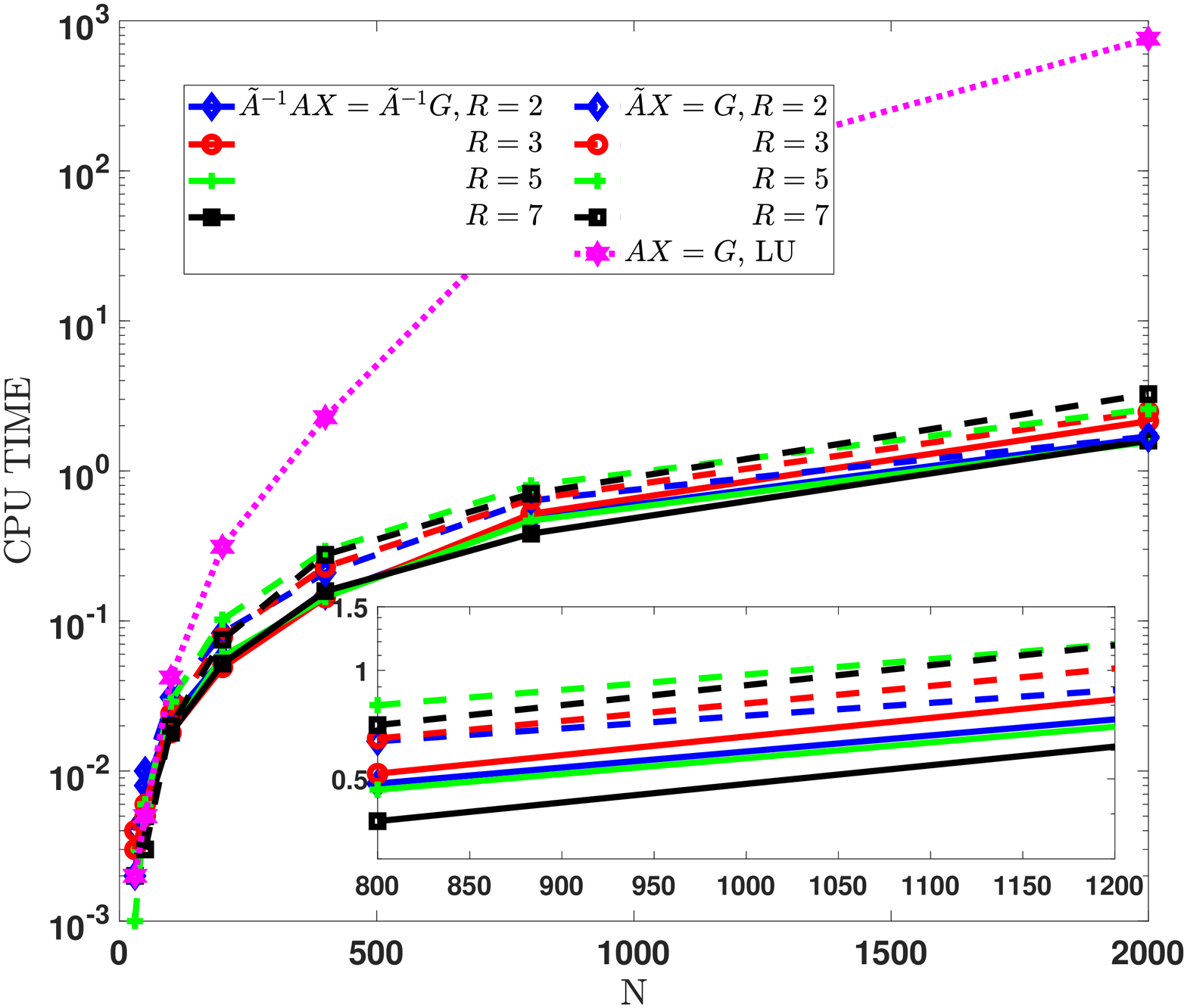}
\caption{\scriptsize Example \ref{ex:uex}:
Accuracy and cost for $h$-refinement. $\rho=100,\,\theta = 1,\, P=3$.  Left: $L^{\infty}$-errors versus number of elements $N$. The preconditioned system has the same accuracy as the original system whereas the H-matrix approximation diverges although initially seems to converge for $R\ge 2$ (see inset).  Right: CPU time versus number of elements $N$. The cases on the left column of the legend are represented by solid lines while the cases on the right column are represented by dash lines. The inset shows a zoom-in plot for clarity.
}
\label{fig5} 
\end{figure}

Similarly, we show the results of $h$-refinement and $p$-refinement for the case $\rho=100$ in Figure \ref{fig5} and \ref{fig6}, respectively. We can see that all results are similar with the ones for the case $\rho= 0$. We then can have the similar conclusion.
We note that the efficient PGS method~\cite{MaoKar18} would not work well in this case because of the presence of the reaction term.

\begin{figure}[!t]
\centering
\includegraphics[width=0.48\textwidth]{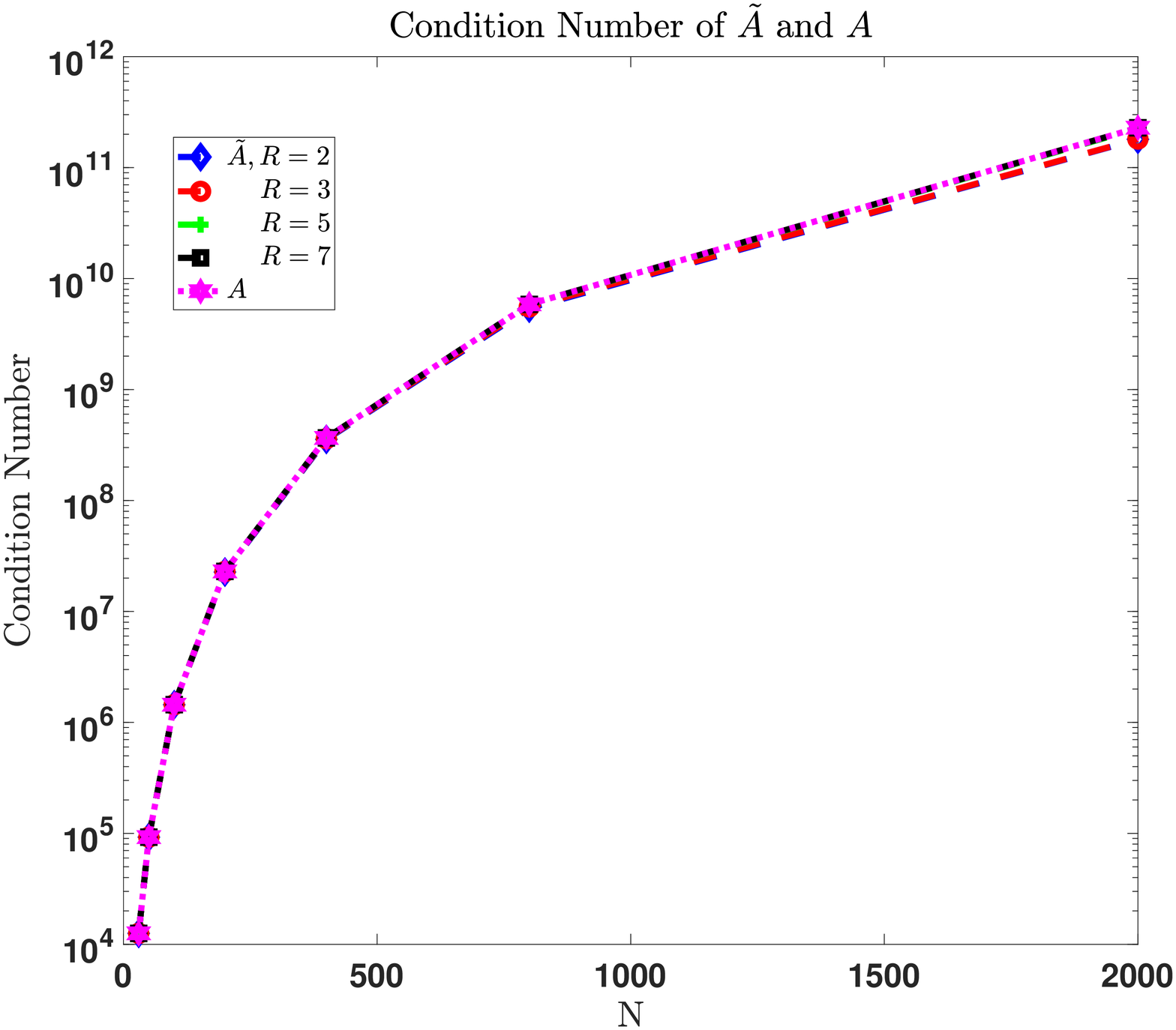}
\includegraphics[width=0.48\textwidth]{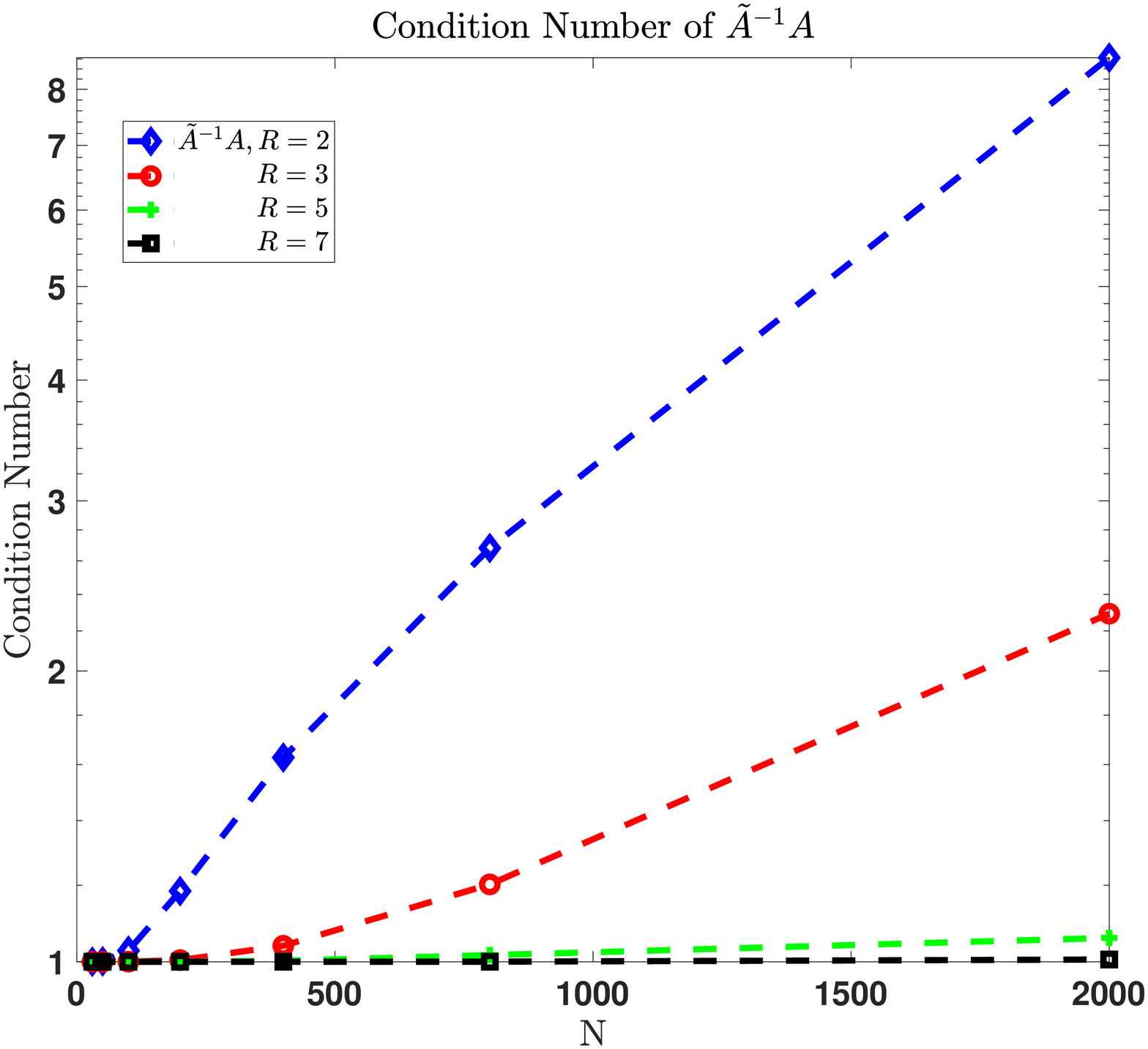}
\caption{\scriptsize Example \ref{ex:uex}:  Condition number for $h$-refinement for different values of rank $R$. $\rho=100,\,\theta = 1,\, P=3$. Left: condition number of $\tilde{A}$ and $A$; right: condition number of $\tilde{A}^{-1}A$.
}
\label{fig52} 
\end{figure}

\begin{figure}[!t]
\centering
\includegraphics[width=0.48\textwidth]{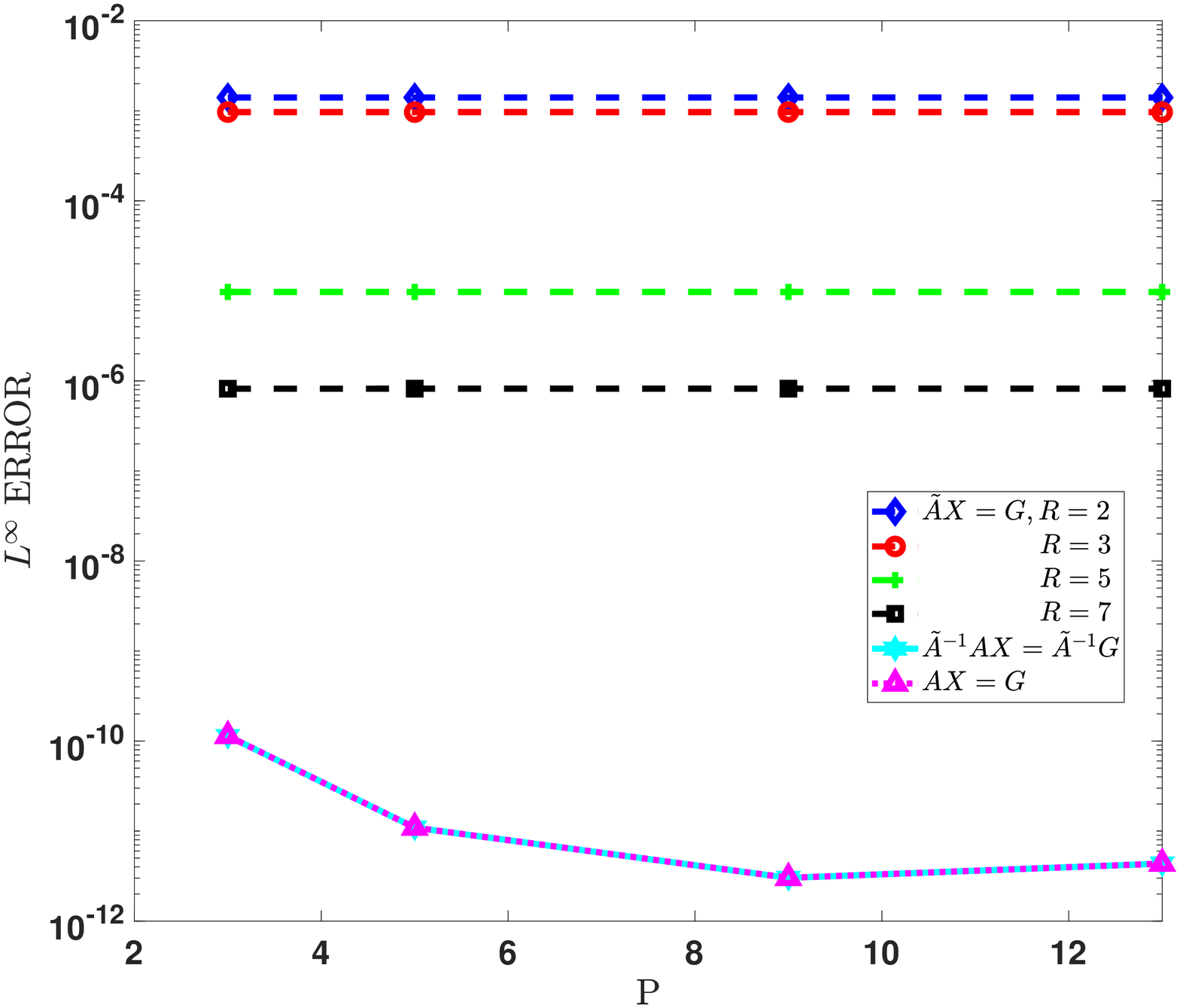}
\includegraphics[width=0.48\textwidth]{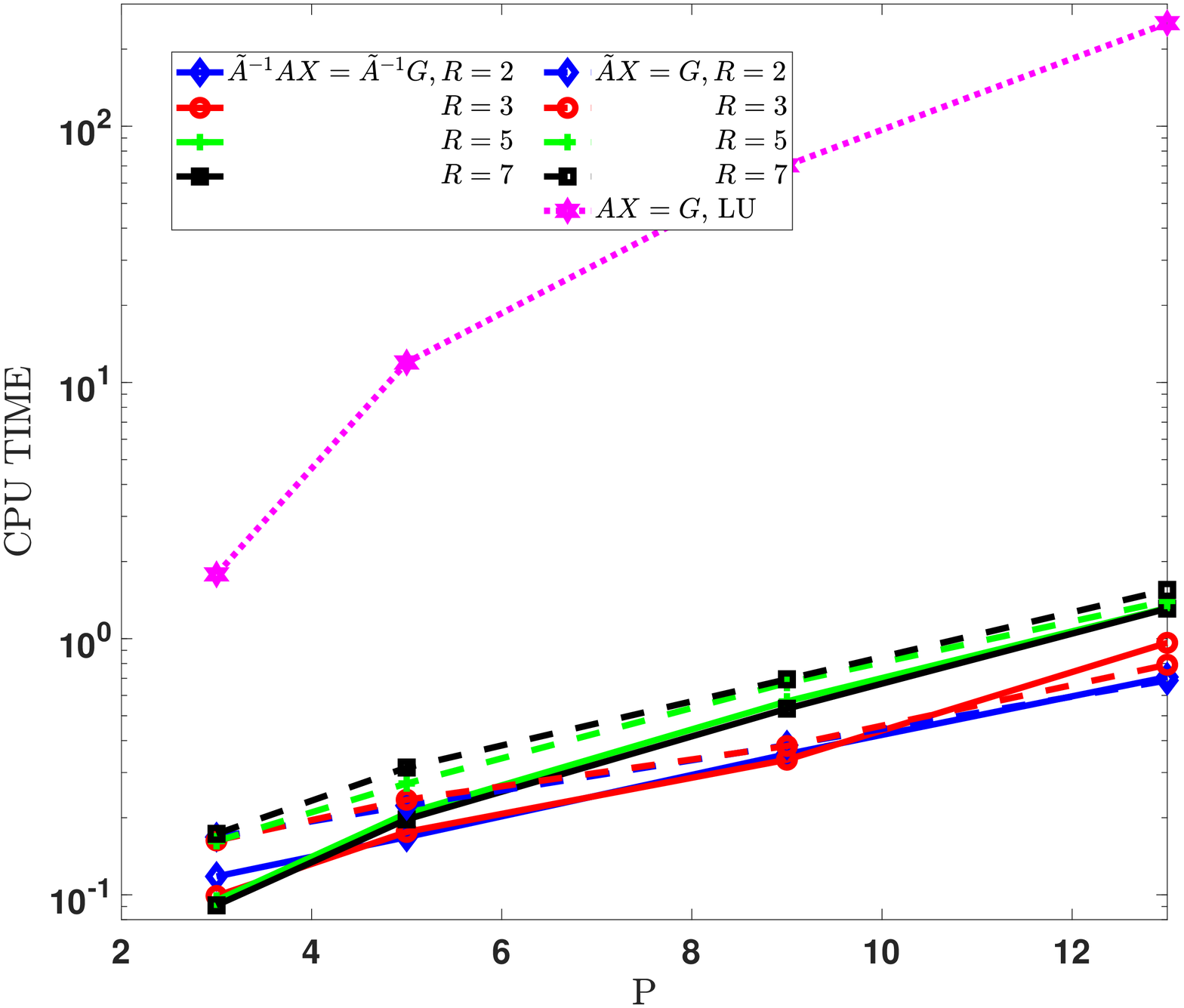}
\caption{\scriptsize Example \ref{ex:uex}:  Accuracy and cost for $p$-refinement. $\rho=100,\,\theta = 1,\, N=300$.  Left: $L^{\infty}$-errors versus polynomial degree $P$. The preconditioned system has the same accuracy as the original system whereas the H-matrix approximation  does not improve with $p$-refinement.  Right: CPU time versus polynomial degree $P$. The cases on the left column of the legend are represented by solid lines while the cases on the right column are represented by dash lines. The inset shows a zoom-in plot for clarity.
}\label{fig6} %
\end{figure}

\begin{figure}[!t]
\centering
\includegraphics[width=0.48\textwidth]{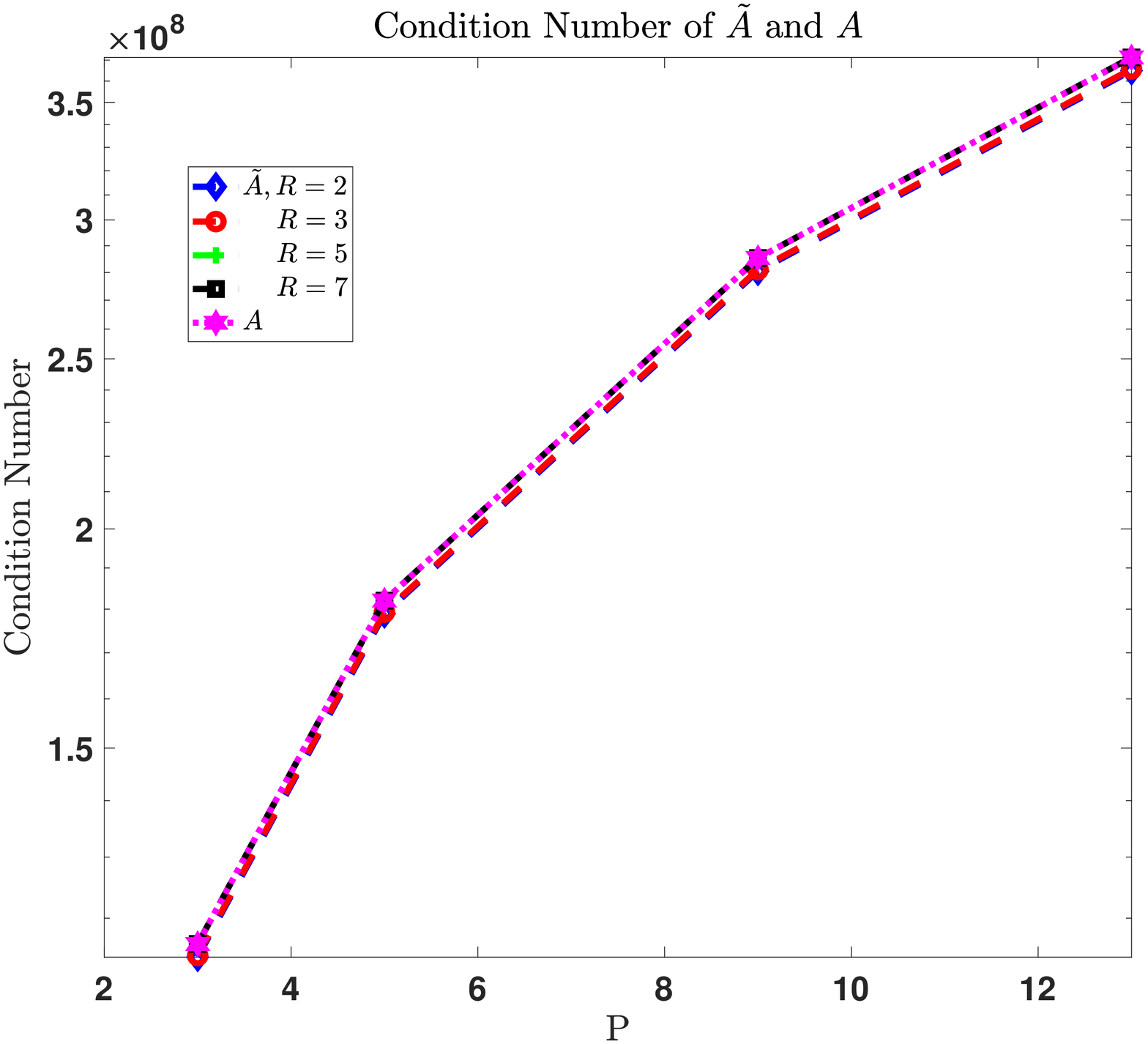}
\includegraphics[width=0.48\textwidth]{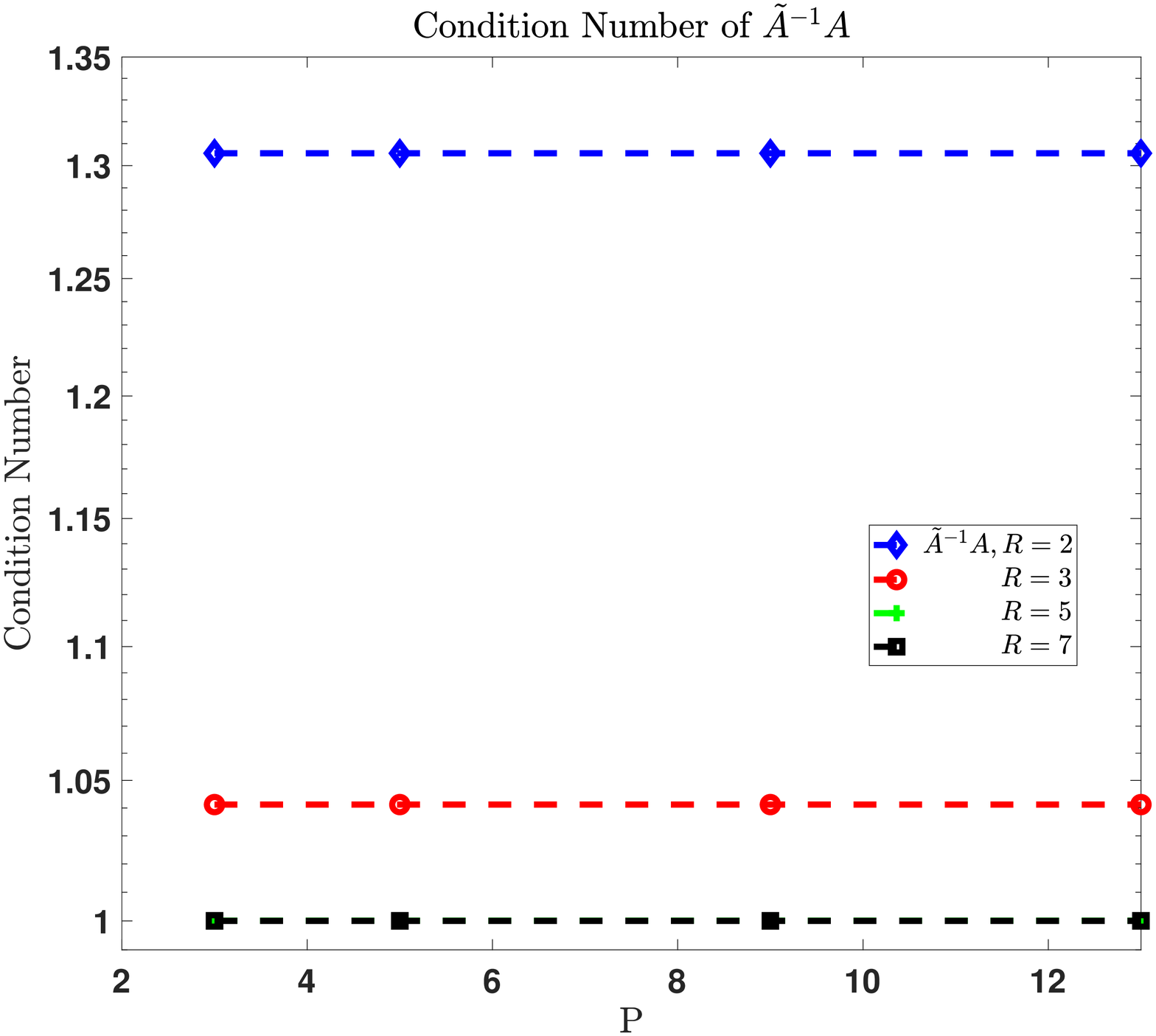}
\caption{\scriptsize Example \ref{ex:uex}:
Condition number for $p$-refinement for different values of rank $R$. $\rho=100,\,\theta = 1,\, N=300$. Left: condition number of $\tilde{A}$ and $A$; right: condition number of $\tilde{A}^{-1}A$.
}\label{fig62} %
\end{figure}

\begin{figure}[!t]
\centering
\includegraphics[width=0.48\textwidth]{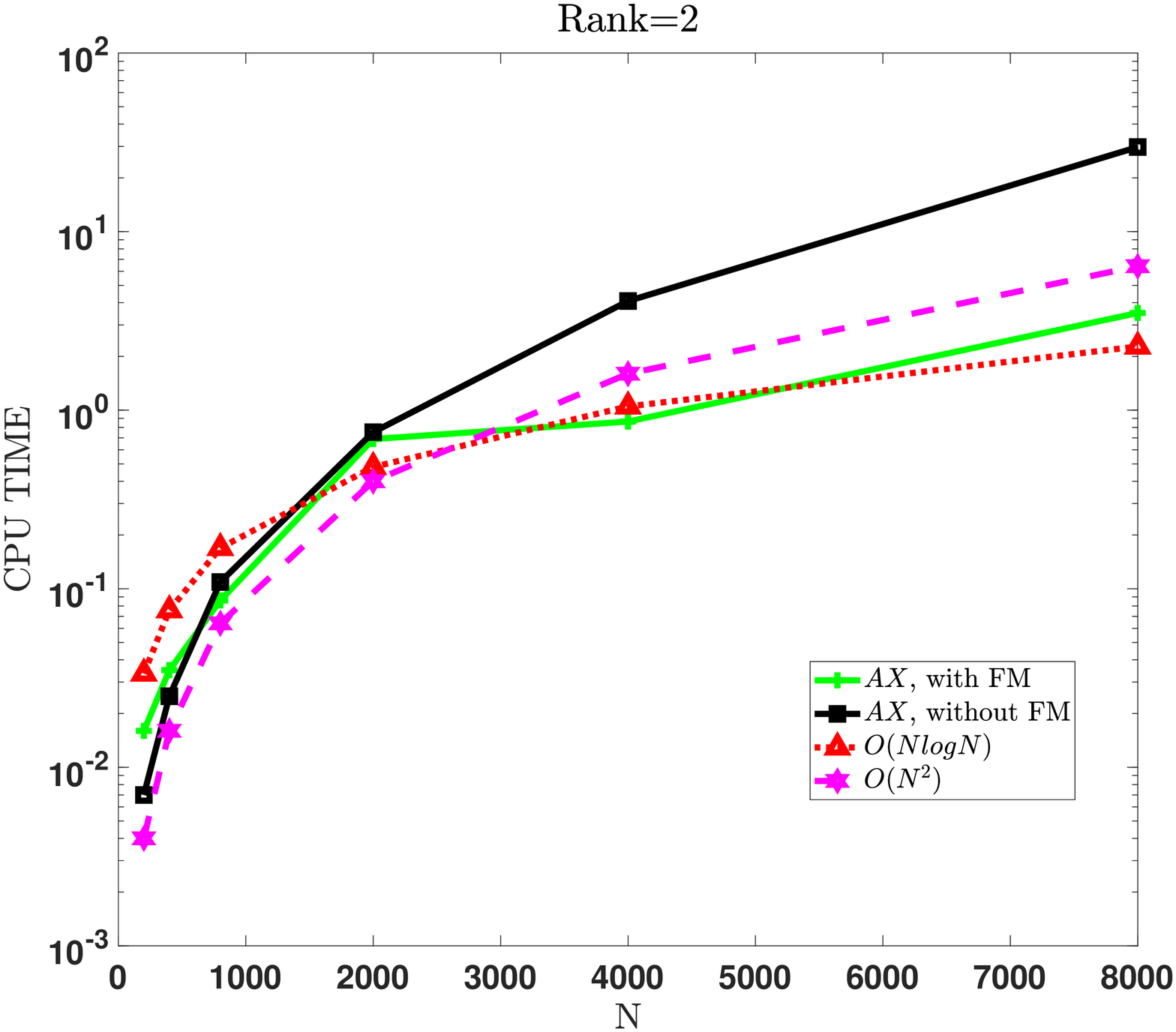}
\includegraphics[width=0.48\textwidth]{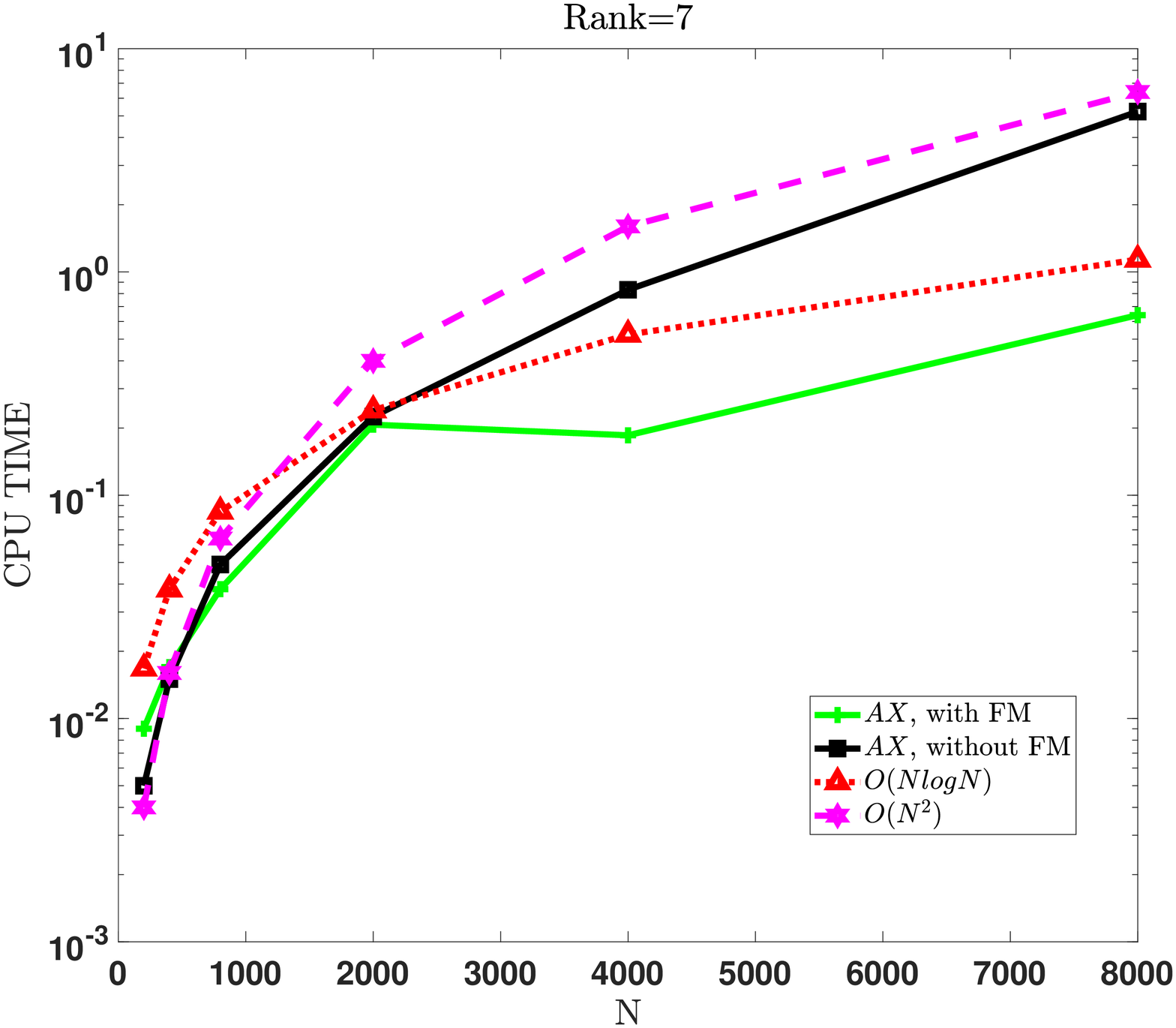}\\
\caption{\scriptsize Example \ref{ex:uex}:  Comparison of CPU time with or without fast matrix-vector multiplication for uniform mesh. $\rho=0,\,\theta = 1,\, P=3$. `FM' denotes fast multiplication. Left: $R = 2$;  right: $R = 7$.
}
\label{fig7} %
\end{figure}

As presented in subsection \ref{sec:toeplitz}, the submatrices of the stiffness matrix $S$ are Toeplitz matrices for a uniform mesh or a geometric mesh. This suggests to use the fast matrix-vector multiplication.
We illustrate this by implementing a numerical simulation using uniform mesh and examine the computational cost of the preconditioned system, where the matrix-vector multiplication $AX$ is computed with or without the fast matrix-vector multiplication. The results are shown in Figure \ref{fig7} for different values of rank $R = 2,7$.
Observe that the computational cost is of order $O(N^2)$ if not using the fast matrix-vector multiplication while it is of order $O(N\log N)$ if using the fast matrix-vector multiplication.
This is in agreement with our previous analysis.

We now consider a example with two-sided fractional derivative and smooth-right-hand function $f(x)$.
\begin{exam}\label{ex:smf}
Let $f(x) = 1$, $\theta = 0.75, \, \rho = 0$,  $a=0,\, b=10$ and $c_1 = c_2 = 0$.
In this case, the exact solution is obtained by a mapping and the equation (44) of \cite{MaoKar18}.
\end{exam}

\begin{figure}[!t]
\centering
\includegraphics[width=0.48\textwidth,height=0.4\textwidth]{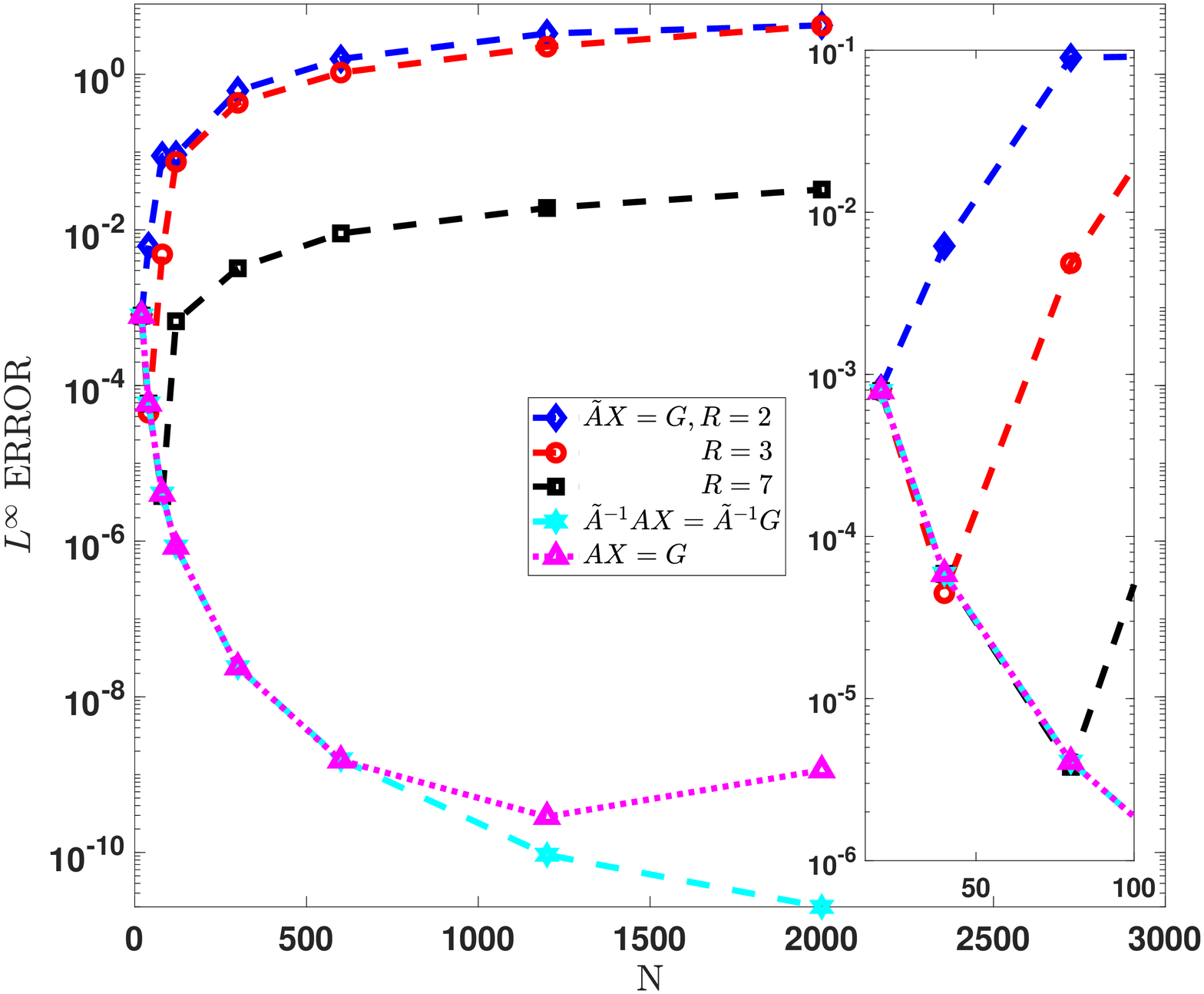}
\includegraphics[width=0.48\textwidth,height=0.4\textwidth]{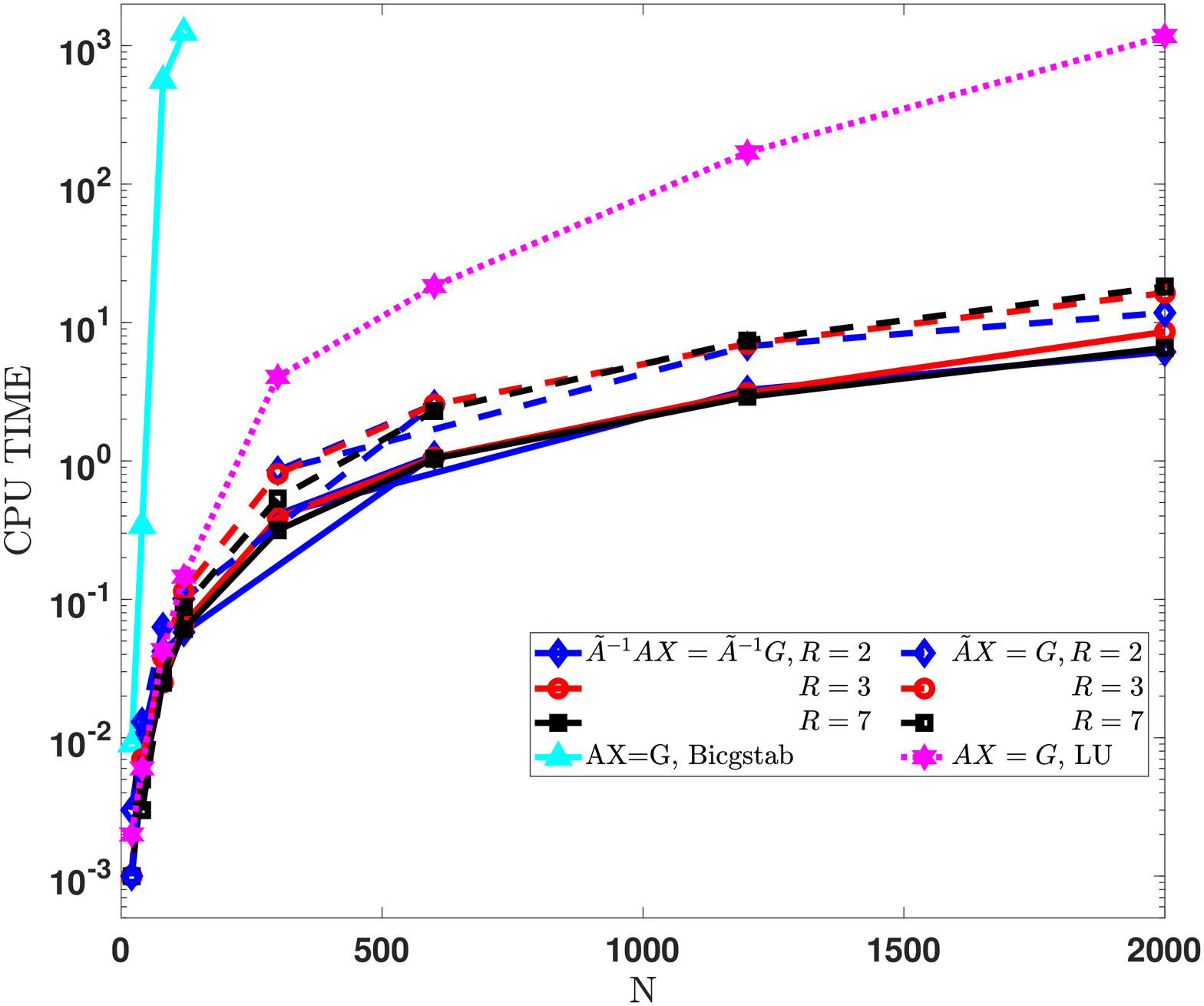}
\caption{\scriptsize Example \ref{ex:smf}:  Accuracy and cost for $h$-refinement. $\rho=0,\,\theta = 0.75,\, P=3$.  Left: $L^{\infty}$-errors versus number of elements $N$. The preconditioned system has the same accuracy as the original system whereas the H-matrix approximation diverges although initially seems to converge for $R\ge 2$ (see inset).  Right: CPU time versus number of elements $N$. The cases on the left column of the legend are represented by solid lines while the cases on the right column are represented by dash lines.
}
\label{figts10} 
\end{figure}

\begin{figure}[!t]
\centering
\includegraphics[width=0.48\textwidth,height=0.4\textwidth]{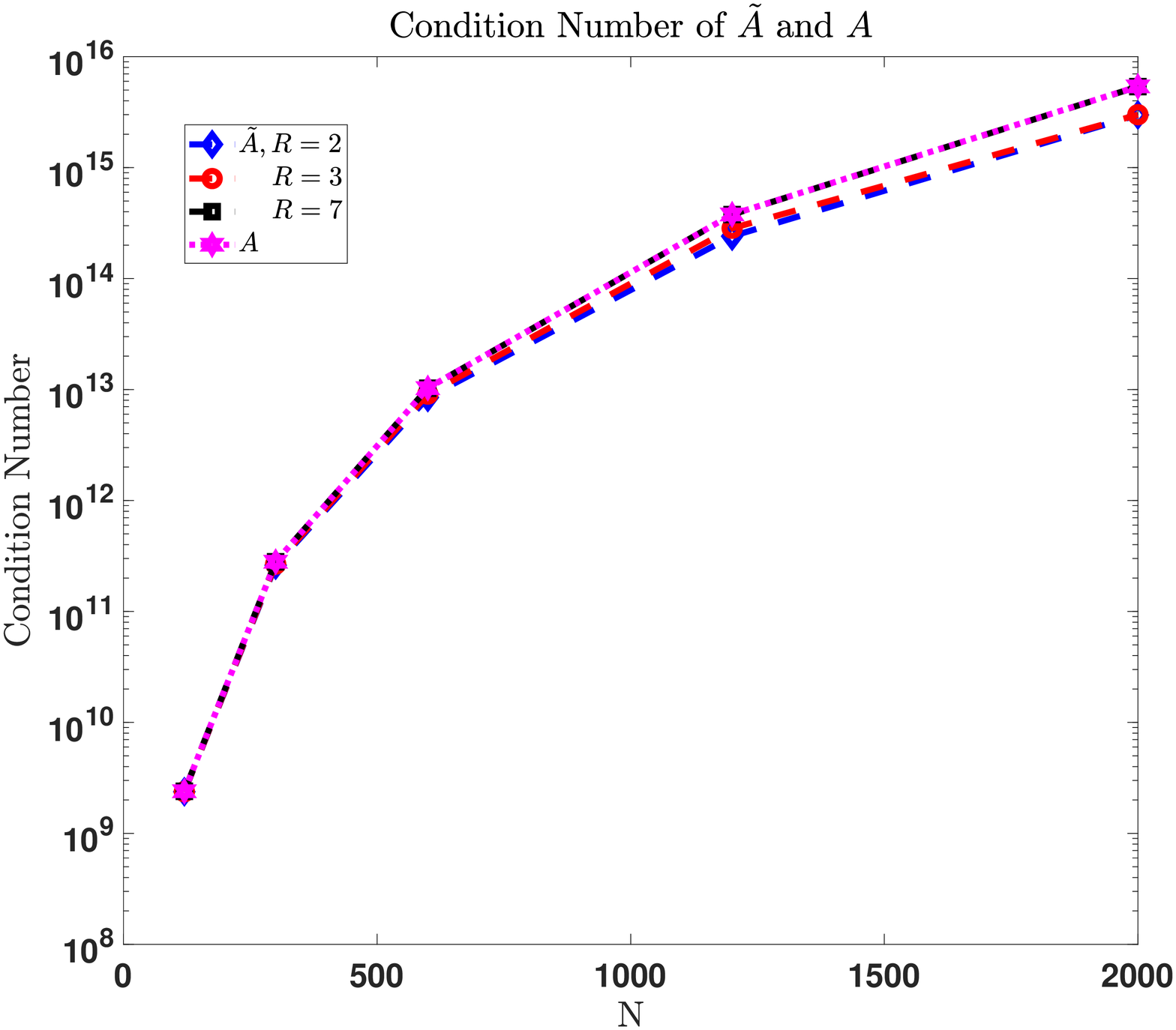}
\includegraphics[width=0.48\textwidth,height=0.4\textwidth]{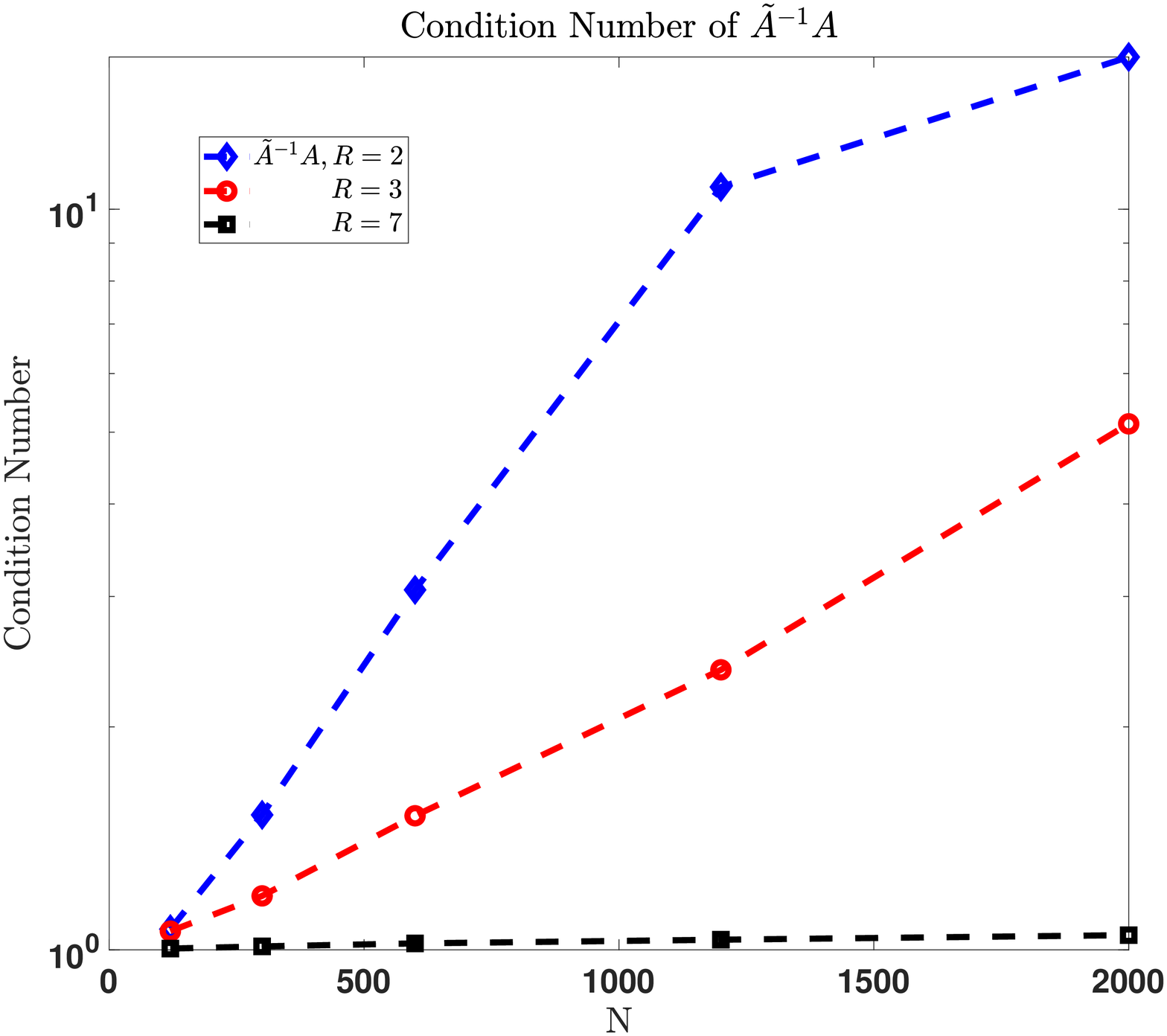}
\caption{\scriptsize Example \ref{ex:smf}:  Condition number for $h$-refinement for different values of rank $R$. $\rho=0,\,\theta = 0.75,\, P=3$. Left: condition numbers of $\tilde{A}$ and $A$; right: condition numbers of $\tilde{A}^{-1}A$.
}
\label{figts20} 
\end{figure}

\begin{figure}[!t]
\centering
\includegraphics[width=0.48\textwidth,height=0.4\textwidth]{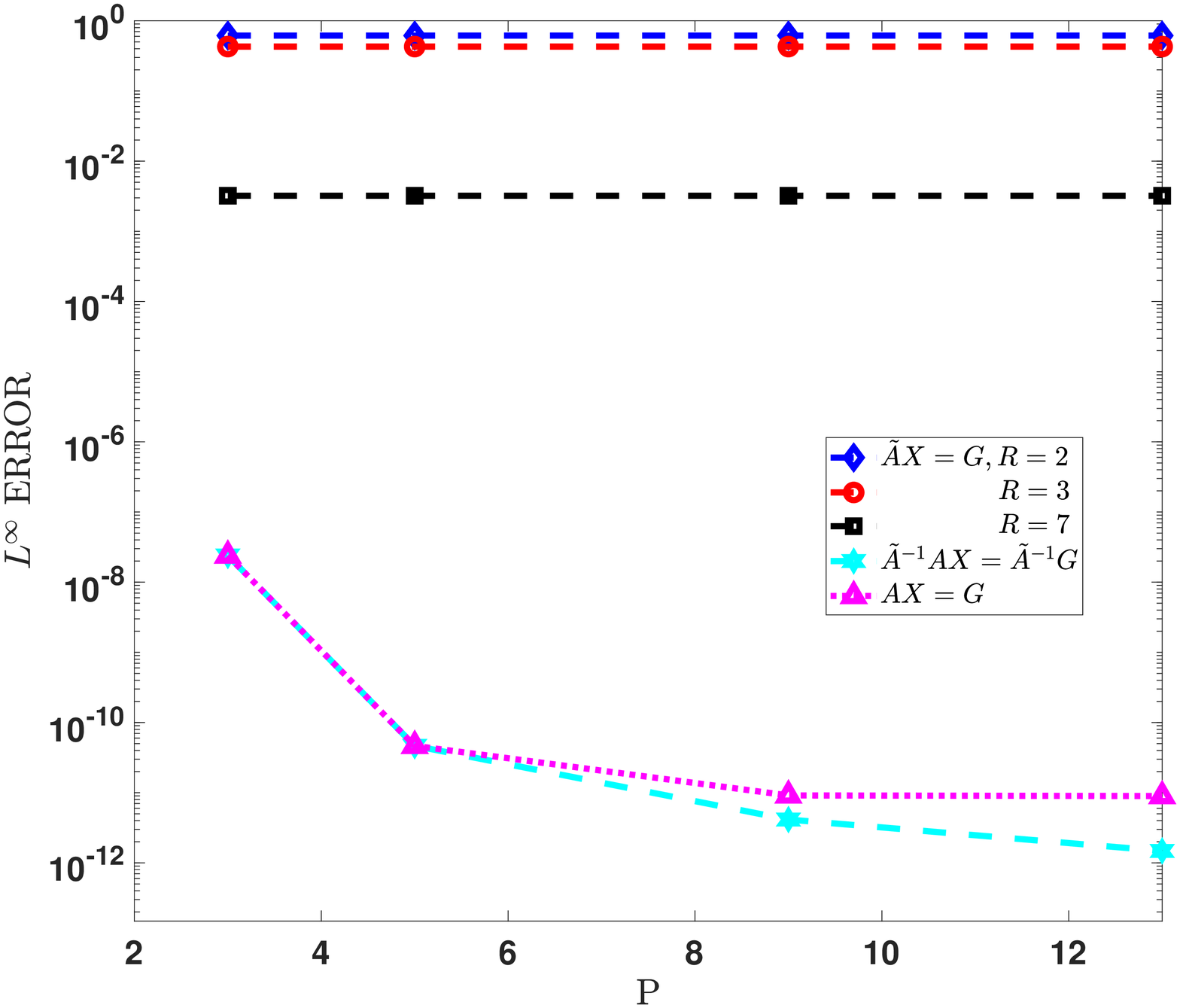}
\includegraphics[width=0.48\textwidth,height=0.4\textwidth]{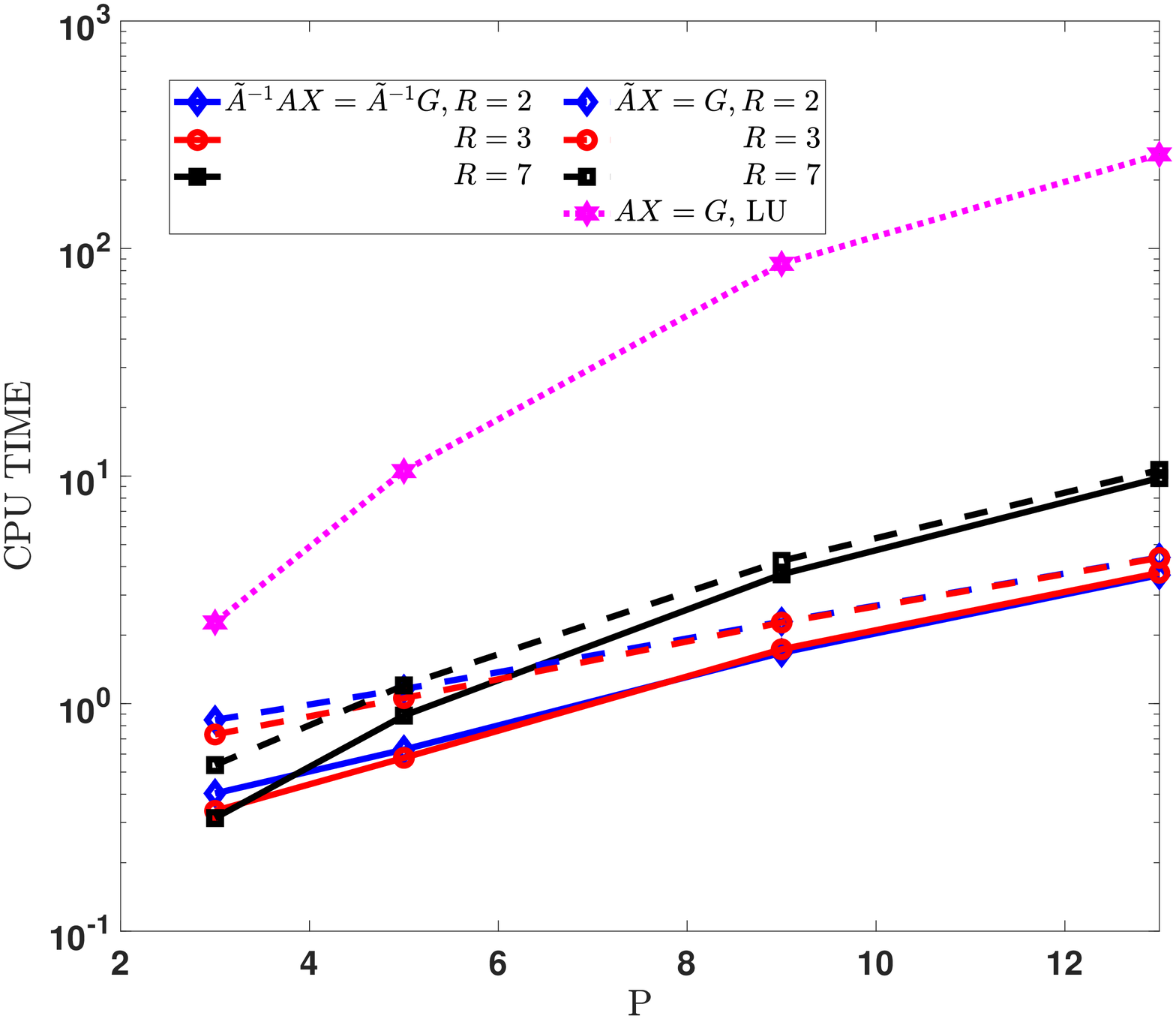}
\caption{\scriptsize Example \ref{ex:smf}: Accuracy and cost for $p$-refinement. $\rho=0,\,\theta = 0.75,\, N=300$.  Left: $L^{\infty}$-errors versus degrees of polynomial $P$. The preconditioned system has the same accuracy as the original system whereas the H-matrix approximation does not improve with $p$-refinement. Right: CPU time versus degrees of polynomial $P$. The cases on the left column of the legend are represented by solid lines while the cases on the right column are represented by dash lines.
}
\label{figts10p} 
\end{figure}

\begin{figure}[!t]
\centering
\includegraphics[width=0.48\textwidth,height=0.4\textwidth]{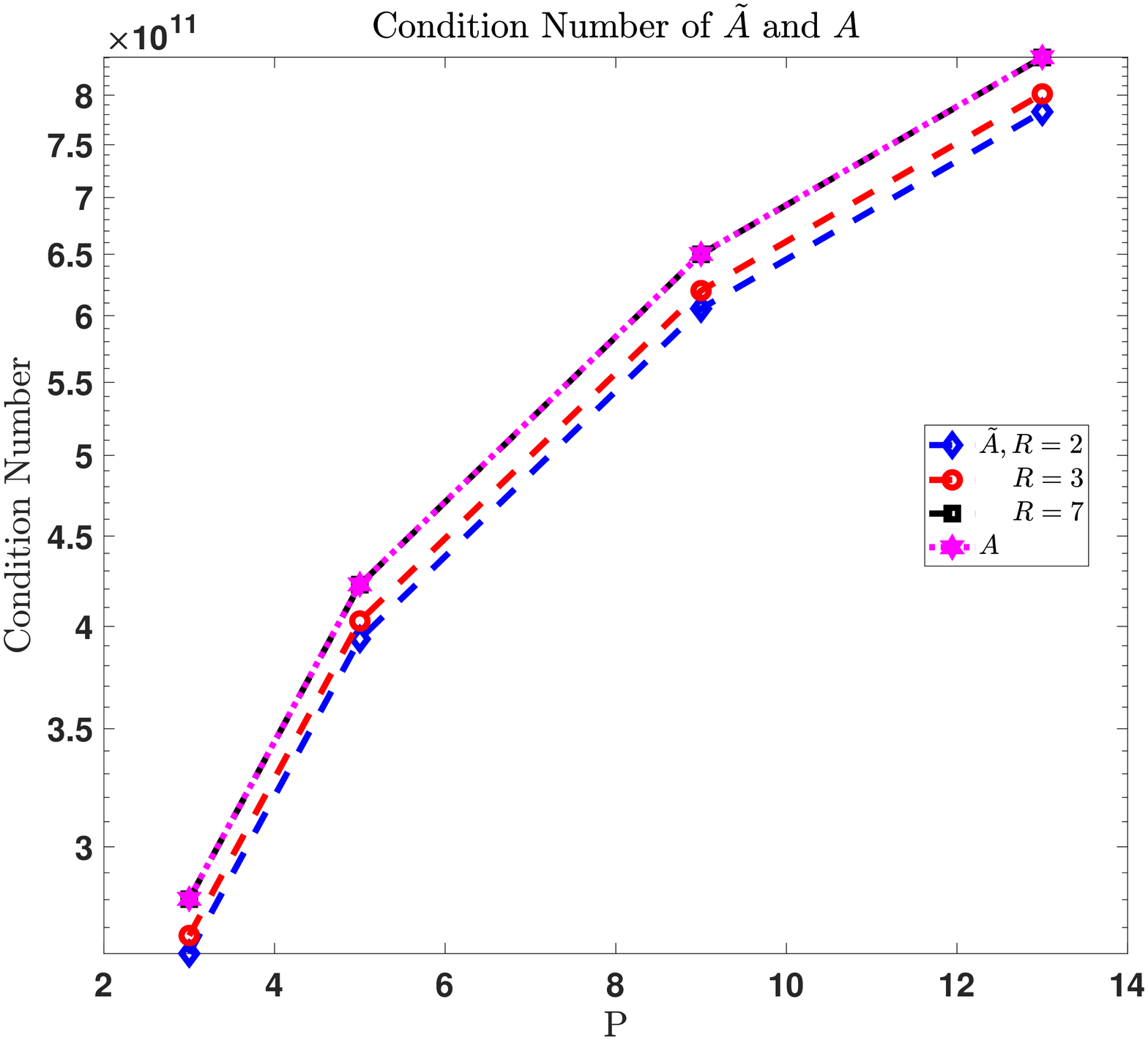}
\includegraphics[width=0.48\textwidth,height=0.4\textwidth]{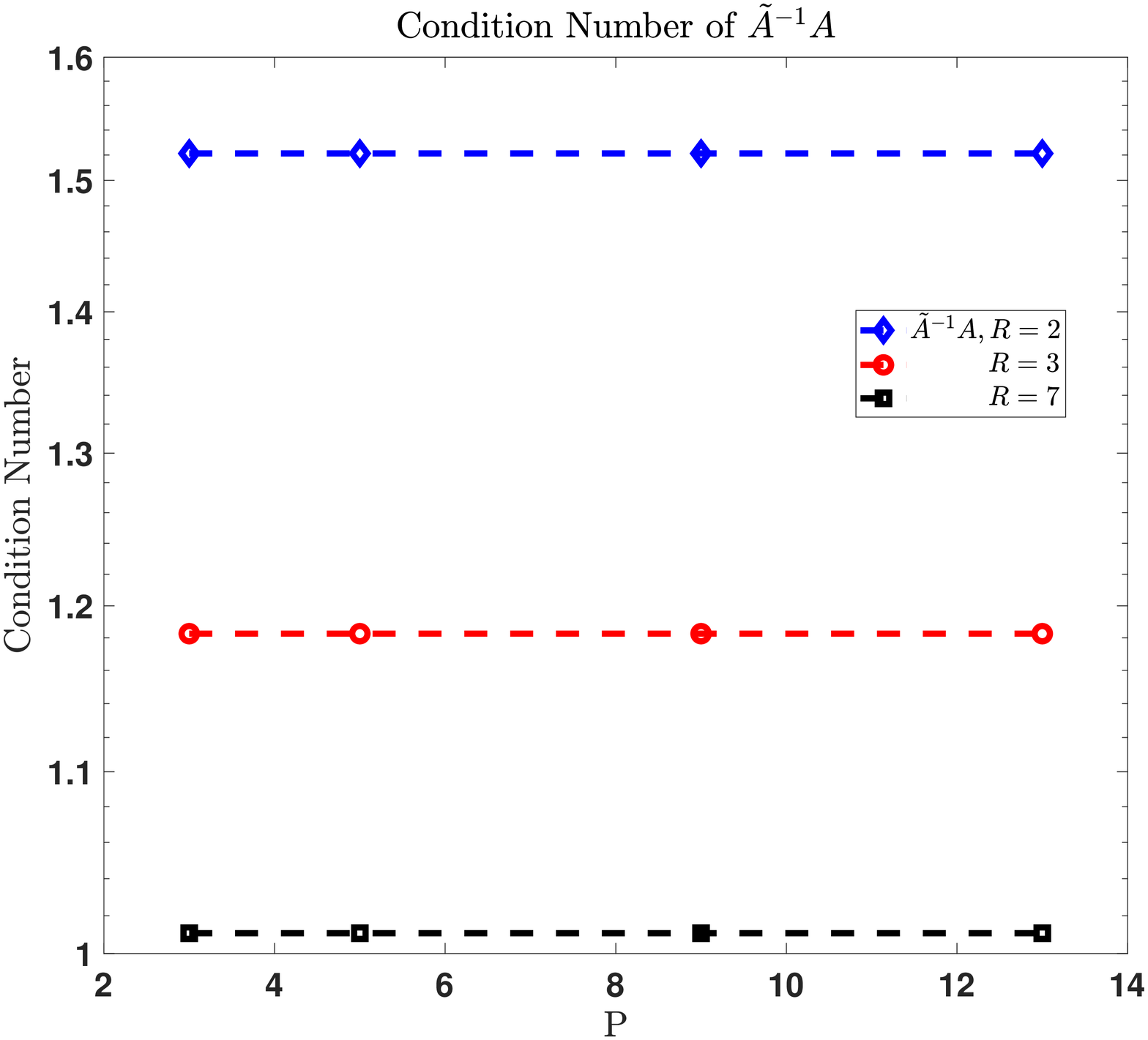}
\caption{\scriptsize Example \ref{ex:smf}:  Condition number for $p$-refinement for different values of rank $R$. $\rho=0,\,\theta = 0.75,\, N=300$. Left: condition number of $\tilde{A}$ and $A$; right: condition number of $\tilde{A}^{-1}A$.
}
\label{figts20p} 
\end{figure}

The mesh used here is the graded mesh, namely, $x_i = 5*(2i/N)^7,  i= 0,1,\ldots, N/2$, $x_{i} = 10-5*((N-2i)/N)^5,i= N/2+1,\ldots,N$. For the $h$-refinement, the results of accuracy and computational cost  are shown in Figure \ref{figts10} while the results of the condition numbers of $A$, $\tilde{A}$ and $\tilde{A}^{-1}A$ are shown in Figure \ref{figts20}.  For the $p$-refinement, the results are shown in Figure \ref{figts10p} (accuracy and computational cost) and Figure \ref{figts20p} (condition numbers), respectively.
Observe that all the results behave similarly as the ones of the previous example for the one-sided problem. This indicates that our algorithm also works for the two-sided problem.

\section{Conclusion}\label{sec:conclusion}
In this work, based on the hierarchical matrix (H-matrix) approximation, we developed a fast solver for the SEM applied to the two-sided fractional diffusion equation.
We first constructed a H-matrix to approximate the stiffness matrix of the SEM by replacing the singular kernel by a degenerate kernel. We also derived the corresponding error estimates.
We then solved efficiently the H-matrix approximation problem with a hierarchical LU decomposition method by which we reduced the computational cost to $O(R^2 N_d  \log^2N) +O(R^3 N_d  \log N)$.
However, we could not obtain high accuracy for the SEM discretization. Thus, to recover the high accuracy of the SEM, we used the H-matrix approximation problem as a preconditioner for the original problem. Numerical results show that the condition number of the preconditioned system is independent of the polynomial degree $P$ and grows with the number of elements, but at modest values of the rank $R$ it is below order 10 in our experiments, which represents a reduction of more than 11 orders of magnitude from the unpreconditioned system. The corresponding cost is $O(R^2 N_d \log^2 N)+O(R^3 N_d  \log N)+O(N_d^2)$. Moreover, by using a uniform mesh, we further reduced the computational cost to  $O(R^2 N_d\log^2 N) +O(R^3 N_d  \log N)+ O(P^2 N\log N)$ for the preconditioned system.
Overall, the preconditioned system possesses both the advantages of the SEM approximation and the H-matrix approximation. In particular, the preconditioned system can be solved efficiently while retaining the high accuracy of the SEM approximation.

\appendix
\section{Matrix form of $C,\,E,\, F$ with a uniform or geometric mesh}\label{sec:apd:mf}
The submatrix $C$  is  a $(N-1)\times (N-1)$  matrix  can be decomposed as
\begin{equation*}
\begin{aligned}
  C = &  \mathrm{diag}(D_{h,\alpha}(1:N-1))\, T_1 \, \mathrm{diag}(D_{h}(1:N-1))\\
      & + \mathrm{diag}(D_{h,\alpha}(2:N))\, T_2 \, \mathrm{diag}(D_{h}(2:N))\\
      & + \mathrm{diag}(D_{h,\alpha}(2:N))\, T_3 \, \mathrm{diag}(D_{h}(1:N-1))\\
      & + \mathrm{diag}([0;D_{h,\alpha}(2:N-1)])\, T_4 \, \mathrm{diag}([D_{h}(2:N-1);0]),
\end{aligned}
\end{equation*}
where $D_h$ and $D_{h,\alpha}$ are given in \eqref{eqn:Dh},
\begin{align*}
  T_1 = &\left[
\ba{cccc}
b_0&&&\\
b_1&b_0&&\\
\vdots&&\ddots&\\
b_{N-2}&\cdots &b_1& b_0
\ea
\right],\quad
T_2 = \left[
\ba{cccc}
e_0&&&\\
e_1&e_0&&\\
\vdots&&\ddots&\\
e_{N-2}&\cdots &e_1& e_0
\ea
\right],
\\
T_3 = &\left[
\ba{cccc}
f_1&f_0&  &\\
f_2&f_1&f_0&\\
\vdots&\ddots&\ddots&\ddots\\
f_{N-1}&\cdots &f_2& f_1
\ea
\right],
\quad
T_4 = \left[
\ba{cccc}
0&&&\\
g_0&0&&\\
\vdots&&\ddots&\\
g_{N-3}&\cdots &g_0& 0
\ea
\right].
\end{align*}
%
with
\begin{align*}
b_0=&\gamma_0 \int_{-1}^1 (\hat{\phi}_L)' \d x \frac{d}{d x} \int_{-1}^x  (x- t)^{1-\alpha}  \hat{\phi}_L \d t =\frac{\gamma_02^{3-\alpha}}{4(3-\alpha)}, \\
b_k=&\gamma_0 \int_{-1}^1 (\hat{\phi}_L)' \d x \frac{d}{d x} \int_{-1}^1  (x-\hat{q}^k t+ 1+2\hat{q}+\cdots+2\hat{q}^{k-1}+\hat{q}^k)^{1-\alpha}  \hat{\phi}_L \d t \\
=&
\frac{\gamma_02^{3-\alpha}}{4}\left[\frac{\zeta_k^{2-\alpha}-\hat{q}^{2-\alpha}\zeta_{k-1}^{2-\alpha}}{-(2-\alpha)\hat{q}^k}+\frac{-\zeta_{k+1}^{3-\alpha}+\zeta_k^{3-\alpha}
+\hat{q}^{3-\alpha}\zeta_{k-1}^{3-\alpha} -\hat{q}^{3-\alpha}\zeta_{k}^{3-\alpha}}{(2-\alpha)(3-\alpha)\hat{q}^{2k}}\right]
\\
%
e_0=& \gamma_0\int_{-1}^1 (\hat{\phi}_R)' \d x \frac{d}{d x} \int_{-1}^x  (x- t)^{1-\alpha}  \hat{\phi}_R \d t=-\frac{\gamma_02^{3-\alpha}}{4(3-\alpha)}, \\
e_k=&\gamma_0 \int_{-1}^1 (\hat{\phi}_R)' \d x \frac{d}{d x} \int_{-1}^1  (x-\hat{q}^k t+ 1+2\hat{q}+\cdots+2\hat{q}^{k-1}+\hat{q}^k)^{1-\alpha}  \hat{\phi}_R \d t\\
=&
\frac{\gamma_02^{3-\alpha}}{4}\left[\frac{\zeta_{k+1}^{2-\alpha}-\hat{q}^{2-\alpha}\zeta_{k}^{2-\alpha}}{-(2-\alpha)\hat{q}^k}+\frac{-\zeta_{k+1}^{3-\alpha}+\zeta_k^{3-\alpha}
+\hat{q}^{3-\alpha}\zeta_{k-1}^{3-\alpha} -\hat{q}^{3-\alpha}\zeta_{k}^{3-\alpha}}{(2-\alpha)(3-\alpha)\hat{q}^{2k}}\right], \\
f_0=&\gamma_0 \int_{-1}^1 (\hat{\phi}_R)' \d x \frac{d}{d x} \int_{-1}^x  (x- t)^{1-\alpha}  \hat{\phi}_L \d t =-\frac{\gamma_02^{3-\alpha}}{4(3-\alpha)}, \\
f_k=& \gamma_0\int_{-1}^1 (\hat{\phi}_R)' \d x \frac{d}{d x} \int_{-1}^1  (x-\hat{q}^k t+ 1+2\hat{q}+\cdots+2\hat{q}^{k-1}+\hat{q}^k)^{1-\alpha}  \hat{\phi}_L \d t=-b_k,\\
%
g_0=&\gamma_0 \int_{-1}^1 (\hat{\phi}_L)' \d x \frac{d}{d x} \int_{-1}^x  (x- t)^{1-\alpha}  \hat{\phi}_R \d t=\frac{\gamma_02^{3-\alpha}}{4(3-\alpha)} , \\
g_k=&\gamma_0 \int_{-1}^1 (\hat{\phi}_L)' \d x \frac{d}{d x} \int_{-1}^1  (x-\hat{q}^k t+ 1+2\hat{q}+\cdots+2\hat{q}^{k-1}+\hat{q}^k)^{1-\alpha}  \hat{\phi}_R \d t=-e_k,
\end{align*}
here
$\hat{\phi}_L(x) = \frac{x+1}{2}, \,\hat{\phi}_R(x)=\frac{1-x}{2},\,x\in[-1,1]$ and $\zeta_k=\frac{1-\hat{q}^{k}}{1-\hat{q}},  k=1, 2,\ldots$.

$F$ is  a $P \times 1$ block matrix with  each block $F^{p}$  being a $N \times (N-1)$ matrix.
Let
\be
\bar{F}^{p}=[F^p,0],
\ee
then $\bar{F}^{p}$ is a $N\times N$ matrix and can be decomposed as
\begin{equation*}
\begin{aligned}
  \bar{F}^{p} = &  \mathrm{diag}(D_{h,\alpha})\, T_5 \, \mathrm{diag}([D_{h}(1:N-1);0])\\
      & + \mathrm{diag}([0;D_{h,\alpha}(2:N)])\, T_6 \, \mathrm{diag}([D_{h}(2:N);0]),
\end{aligned}
\end{equation*}
where
\begin{equation*}
  T_5 =  \left[
\ba{cccc}
c^{p}_0&&&\\
c^{p}_1&c^{p}_0&&\\
\vdots&&\ddots&\\
c^{p}_{N-1}&\cdots &c^{p}_1& c^{p}_0
\ea
\right],
\quad
T_6 = \left[
\ba{cccc}
0&&&\\
d^{p}_0&0&&\\
\vdots&&\ddots&\\
d^{p}_{N-2}&\cdots &d^{p}_0& 0
\ea
\right]
\end{equation*}
with
%
%
\begin{align*}
c^{p}_0=&\gamma_0 \int_{-1}^1 (\psi_p)' \d x \frac{d}{d x} \int_{-1}^x  (x- t)^{1-\alpha}  \hat{\phi}_L \d t  , \\
c^{p}_k=&\gamma_0 \int_{-1}^1 (\psi_p)' \d x \frac{d}{d x} \int_{-1}^1 (x-\hat{q}^k t+ 1+2\hat{q}+\cdots+2\hat{q}^{k-1}+\hat{q}^k)^{1-\alpha} \hat{\phi}_L \d t,\\
d^{p}_0=&\gamma_0 \int_{-1}^1 (\psi_p)' \d x \frac{d}{d x} \int_{-1}^x  (x- t)^{1-\alpha}  \hat{\phi}_R \d t  , \\
d^{p}_k=&\gamma_0 \int_{-1}^1 (\psi_p)' \d x \frac{d}{d x} \int_{-1}^1 (x-\hat{q}^k t+ 1+2\hat{q}+\cdots+2\hat{q}^{k-1}+\hat{q}^k)^{1-\alpha} \hat{\phi}_R \d t,\\
&k=1, 2,\ldots.
\end{align*}

$E$ is  a $1\times P$ block matrix with each block $E^{q}$ being a $(N-1) \times N$ matrix.
Let
\be
\bar{E}^{q}=
\left[
\begin{array}{c}
0\\
E^q
\end{array}
\right],
\ee
then $\bar{E}^{q}$ is a $N\times N$ matrix and can be decomposed as
then $\bar{F}^{q}$ is a $N\times N$ matrix and can be decomposed as
\begin{equation*}
\begin{aligned}
  \bar{E}^{q} = &  \mathrm{diag}([0;D_{h,\alpha}(2:N)])\, T_7 \, \mathrm{diag}(D_{h})\\
      & + \mathrm{diag}([0;D_{h,\alpha}(1:N-1)])\, T_8 \, \mathrm{diag}([D_{h}(1:N-1);0]),
\end{aligned}
\end{equation*}
where
\begin{equation*}
  T_7 =  \left[
\ba{cccc}
l^{q}_0&&&\\
l^{q}_1&c^{q}_0&&\\
\vdots&&\ddots&\\
l^{q}_{N-1}&\cdots &l^{q}_1& l^{q}_0
\ea
\right],
\quad
T_8 = \left[
\ba{cccc}
0&&&\\
s^{q}_0&0&&\\
\vdots&&\ddots&\\
s^{q}_{N-2}&\cdots &s^{q}_0& 0
\ea
\right]
\end{equation*}
with
\begin{align*}
l^{q}_0=&\gamma_0 \int_{-1}^1 \hat{\phi}_R\d x \frac{d}{d x} \int_{-1}^x  (x- t)^{1-\alpha}(\psi_q)'   \d t , \\
l^{q}_k=&\gamma_0 \int_{-1}^1 \hat{\phi}_R\d x \frac{d}{d x} \int_{-1}^1 (x-\hat{q}^k t+ 1+2\hat{q}+\cdots+2\hat{q}^{k-1}+\hat{q}^k)^{1-\alpha} (\psi_q)' \d t, \\
s^{q}_0=&\gamma_0 \int_{-1}^1 \hat{\phi}_L\d x \frac{d}{d x} \int_{-1}^x  (x- t)^{1-\alpha} (\psi_q)'   \d t , \\
s^{q}_k=&\gamma_0 \int_{-1}^1 \hat{\phi}_L\d x \frac{d}{d x} \int_{-1}^1  (x-\hat{q}^k t+ 1+2\hat{q}+\cdots+2\hat{q}^{k-1}+\hat{q}^k)^{1-\alpha} (\psi_q)' \d t, \\
&k=1, 2,\ldots.
\end{align*}

\bibliographystyle{siamplain} 

\end{document}